\documentclass[11pt]{amsart}
\usepackage{amsmath,amsxtra,amssymb}

\newcommand{\8}{\infty}

\newcommand{\kla}{\left ( }
\newcommand{\mer}{\right ) }

\renewcommand{\for}{\begin{eqnarray*}}
\newcommand{\mel}{\end{eqnarray*}}

\newcommand{\kl}{\pl \le \pl}

\newcommand{\lel}{\pl = \pl}

\newcommand{\lz}{\vspace{0.3cm}}

\newcommand{\nz}{{\mathbb N}}
\newcommand{\nen}{n \in \nz}

\newcommand{\rz}{{\mathbb R}}

\newcommand{\cz}{{\mathbb C}}

\newcommand{\ten}{\otimes}

\newcommand{\pl}{\hspace{.1cm}}
\newcommand{\pll}{\hspace{.3cm}}

\newcommand{\al}{\alpha}
\newcommand{\si}{\sigma}
\newcommand{\Si}{\Sigma}

\newcommand{\la}{\lambda}
\newcommand{\eps}{\varepsilon}

\newcommand{\F}{{\mathcal F}}
\newcommand{\E}{{\mathcal E}}

\newcommand{\A}{{\mathcal A}}
\newcommand{\B}{{\mathcal B}}

\newcommand{\M}{{\mathcal M}}

\newcommand{\N}{{\mathcal N}}

\newcommand{\noo}{\left \|}
\newcommand{\rrm}{\right \|}

\newcommand{\intt}{\int\limits}
\newcommand{\summ}{\sum\limits}

\newtheorem{lemma}{Lemma}[section]
\newtheorem{prop}[lemma]{Proposition}
\newtheorem{theorem}[lemma]{Theorem}
\newtheorem{cor}[lemma]{Corollary}

\newtheorem{rem}[lemma]{Remark}

\newtheorem{fact}[lemma]{Fact}
\newcommand{\re}{\begin{rem}\rm}
  \newcommand{\mar}{\end{rem}}

\begin{document}

\title{A noncommutative version of the John-Nirenberg
theorem}

\author{Marius Junge$^{(1)}$ \and Magdalena Musat$^{(2)}$}
\address{$^{(1)}$Department of Mathematics, 1409 West Green Street,
University of Illinois at Urbana-Champaign, Urbana, IL 61801.}
\email{junge@math.uiuc.edu}
\address{
$^{(2)}$Department of Mathematics, 0112, 9500 Gilman Drive,
University of California , San Diego, La Jolla, CA 92093-0112.}
\email{mmusat@math.ucsd.edu}

\footnotetext{Junge is partially supported by the National Science
Foundation, DMS-0301116.}

\keywords{noncommutative $L_p$-spaces and $BMO$;\ noncommutative
martingales; \ interpolation} \subjclass[2000]{Primary: 46L52;
Secondary: 60G46}

\date{}
\begin{abstract}
We prove a noncommutative version of the John-Nirenberg theorem
for non-tracial filtrations of von Neumann algebras. As an
application, we obtain an analogue of the classical large
deviation inequality for elements of the associated $BMO$ space.
\end{abstract}

\maketitle

\section{Introduction}

The John-Nirenberg theorem is an important tool in analysis and
probability. It provides a characterization of $BMO$\,, the space
of functions of {\em bounded mean oscillation}. In its first
version it was proved by John and Nirenberg \cite{Jo} in 1961.
Through its connection with the theory of Hardy spaces, the
John-Nirenberg result has many applications in harmonic analysis
and the theory of singular integrals (see, e.g., Stein \cite{Ste1,
Ste2}). We refer to the books of Bennett and Sharpley \cite{BeSh},
Garcia-Cuerva and Rubio de Francia \cite{GR}, Garnett \cite{Gar}
and Koosis \cite{Koo} for the interval (function space) version,
respectively to the work of Bass \cite{Ba}, Garsia \cite{Ga} and
Petersen \cite{Pet} for the martingale version of this theorem.

In this paper we analyze analogues of the John-Nirenberg results
in a noncommutative setting. We first recall the classical
results. Let $(\Omega, {\F}, \mathbb{P})$ be a probability space
and $({\F}_n)_{n\geq 0}$ an increasing sequence of
sub-$\sigma$-algebras of ${\F}$. For $1\leq p< \infty$ consider
the norms
\begin{equation}\label{eq3433}
\|x\|_{{BMO}_p}:=\sup\limits_n\|\mathbb{E}((|x-x_{n-1}|^p)|\F_n)\|^{\frac1p}_\infty\,,
\end{equation}
where $x_n=\mathbb{E}(x|\F_n)$\,, for all non-negative integers
$n$\,. The usual $BMO$ norm corresponds to $p=2$ above, i.e.,
$\|x\|_{BMO}=\|x\|_{{BMO}_2}$\,. The John-Nirenberg theorem yields
universal constants $C_1, C_2>0$ with the following property. If
$\|x\|_{{BMO}_1}< C_2$\,, then
\begin{equation}\label{JN3}
\sup\limits_n\,\,\|{\mathbb
{E}}(e^{C_1|x-x_{n-1}|}|\F_n)\|_\infty< 1\,.
\end{equation}
Using the power series expansion for the exponential function, it
follows that \eqref{JN3} is equivalent to
\begin{equation}\label{JN4}
\|\mathbb{E}(|x-x_{n-1}|^p|\F_n)\|^{\frac1p}_\infty\leq Cp\
\|x\|_{BMO}\,,
\end{equation}
for all $n\geq 0$ and all $1\leq p< \infty\,,$ where $C> 0$ is a
universal constant. Furthermore, a standard duality argument
yields the equality
\begin{equation}\label{JN6}
\|\mathbb{E}(|x-x_{n-1}|^p|\F_n)\|^{\frac1p}_\infty=\sup_{{a\in
L_p({\F}_n)}, \,{\|a\|_p}\leq 1} \|(x-x_{n-1})a\|_p\,.
\end{equation}
It can be easily seen that condition (\ref{JN4}) implies that, for
$1\leq p< \infty$\,,
\begin{equation}\label{JN5}
\|x\|_p\leq Cp\ \|x\|_{BMO}\,.
\end{equation}
Moreover, the John-Nirenbeg theorem follows from the validity of
(\ref{JN5}) for arbitrary probability measures and filtrations.
Indeed, let $1\leq p< \infty$\,, and fix an integer $n\geq 0$\,.
Consider a positive element $a\in L_p(\F_n)$ with $\|a\|_p=1$\,,
and denote by $\mu_a$ the probability measure defined by
$d\mu_a=a^{\frac1p}d\mathbb{P}$\,. Then
\begin{equation}\label{eq:newmeas}
\|(x-x_{n-1})a\|_p=\|x-x_{n-1}\|_{L_p(\mu_a)}\,.
\end{equation}
Denote by $BMO((\F_k)_{k\geq n}\,, \mu_a)$ the $BMO$ space
associated to the triplet $(\Omega, \F, \mu_a)$ and the filtration
$(\F_k)_{k\geq n}$ of $\F$\,. Applying (\ref{JN5}) for the
probability measure $\mu_a$ and the new filtration $(\F_k)_{k\geq
n}$ of $\F$\,, it follows that
\begin{equation}\label{eq:jnewfiltr}
\|x-x_{n-1}\|_{L_p(\mu_a)}\leq  Cp \
\|x-x_{n-1}\|_{BMO((\F_k)_{k\geq n}, \,\mu_a)}\,.
\end{equation}
Note that $\|x-x_{n-1}\|_{BMO((\F_k)_{k\geq n}\,, \mu_a)}\leq
\|x\|_{BMO((\F_k)_{k\geq 0}\,, \mathbb{P})}\,.$ An application of
(\ref{eq:newmeas})\,, together with formula (\ref{JN6}) yields now
the inequality (\ref{JN4})\,.

In the form (\ref{JN5}), one can formulate a noncommutative
version of the John-Nirenberg theorem. First this requires an
appropriate definition of $BMO$ spaces. In their seminal paper
\cite{PX}, Pisier and Xu proved the noncommutative analogues of
the Burkholder-Gundy square function inequalities. Fermionic
versions of the square function inequalities had previously been
considered by Carlen and Kr\'ee \cite{CK}. In \cite{PX}, Pisier
and Xu also introduced the $BMO$ space for noncommutative
martingales and proved the analogue of the classical
Fefferman-Stein duality $BMO=(H_1)^*$. Their work triggered a
rapid development in the $L_p$-theory of noncommutative
martingales, see, e.g., Junge \cite{Ju}, Junge and Xu \cite{JX,
JX2, JX1} and Randrianantoanina \cite{Ra, Ra2}.

In the following we assume that $\N$ is a von Neumann algebra,
$\phi$ a normal faithful state on $\N$ and $(\N_n)_{n\geq 0}$ an
increasing sequence of von Neumann subalgebras whose union
generates $\N$ in the w$^*$-topology.
Moreover, we assume that for all positive integers $n$ there exist
normal conditional expectations $\E_n:\N\to \N_n$ such that
\begin{equation}\label{e.1.45679}
\phi(\E_n(x)y)=\phi(xy)
\end{equation}
for all $x\in \N$ and $y\in \N_n$. We recall the following
definitions (see \cite{JX}, \cite{PX})
 \[ \noo x\rrm_{BMO^c} \lel \sup_m \ \sup_{n\leq m} \noo
 \E_n((x_m-x_{n-1})^*(x_m-x_{n-1}))\rrm_{\infty}^{\frac12} \pl ,\]
where $x_n=\E_n(x)$\,, for all $n\geq 0\,.$ Respectively,
 \[ \noo x\rrm_{BMO}\lel \max\{\noo x\rrm_{BMO^c},\noo
 x^*\rrm_{BMO^c} \} \pl .\]
The conditioned $L_\infty$-spaces, $L_\infty^c(\N, \E_n)$
associated to the conditional expectations $\E_n$ were introduced
in \cite{Ju} (see also \cite{Ri}) by defining
 \[ \noo x\rrm_{L_\infty^c(\N,\E_n)} \lel \noo \E_n(x^*x)\rrm_{\infty}^{\frac12}
 \pl .\]
We shall point out that in general, for $2< p< \infty$\,, the
expression $\|\E_n((x^*x)^{\frac{p}{2}})\|^{\frac1p}_\infty$ does
not necessarily provide a norm in the noncommutative setting.
Therefore, we will use a different approach to generalize the
$BMO_p$-norms. Namely, motivated by (\ref{JN6}), we define for
$2\leq p\leq \infty$\,,
\begin{eqnarray}
\noo x\rrm_{BMO_p^c}&=&\sup_m \ \sup_{n\leq m} \ \sup_{a\in
 L_p(N_n), \,\noo a\rrm_p\le 1} \noo (x_m-x_{n-1})a\rrm_p\,.\label{eq4534}
 \end{eqnarray}
Respectively, we define
\begin{eqnarray}\label{eq78789}
\noo x\rrm_{BMO_p}&\lel& \max\{\noo x\rrm_{BMO_p^c},\noo
 x^*\rrm_{BMO_p^c}\}\,.
\end{eqnarray}
These norms can, in fact, be obtained by interpolation between
conditional $L_p$-spaces. More details are given in Section 4.
(See \cite{JPa} for a more general discussion of such norms.)

Our martingale version of the John-Nirenberg theorem reads as
follows.

\begin{theorem}\label{martingal} There exists a universal constant $c> 0$ such that for all $2< p<\infty$\,,
\begin{equation}\label{JN9}
\noo x\rrm_{BMO}\kl \noo x\rrm_{BMO_p} \kl cp\noo
 x\rrm_{BMO} \pl .
 \end{equation}
\end{theorem}

The order $cp$ of the constant in the inequality (\ref{JN9}) above
is the same as in the commutative setting. In a preliminary
version of this paper we were only able to prove that (\ref{JN9})
holds with a constant of the order $cp^2$\,. The right order $cp$
was obtained using very recent results of Randrianantoanina
\cite{Ra2} on the optimal order of growth for the constants in the
noncommutative Burkholder-Gundy square function inequalities. For
the proof of Theorem \ref{martingal}, note that formula
(\ref{eq4534}), which is the key point in the change of density
argument described in classical setting, leads to the passage from
traces to various states on the von Neumann algebra. Therefore we
consider martingales with respect to states and their modular
theory.

As an application of a noncommutative version of Chebychev
inequality proved by Defant and Junge \cite{DJ}, we obtain an
analogue of the classical large deviation inequality for elements
of $BMO$. Namely, if $\ \|x\|_{BMO}< 1$\,, then
\begin{equation}\label{eqlargedi}
\mathbb{P}\{|x-x_0|> t\}< {C_2}e^{-{t} {C_1}}\,.
\end{equation}
for all $t> 0$\,, where $C_1, C_2$ are universal constants (see
\cite{Gar}, \cite {Ga}). More precisely, we prove
\begin{theorem}\label{noncomldiv}
There exist universal constants $c_1, c_2> 0$ such that if
$\|x\|_{BMO}< 1$\,, then for
all $t> 0\,,$ there exists a projection $f\in \N$ such that
$\|(x-x_0)f\|\leq t$ and
\begin{equation}\label{eq:ldiv}
\phi(1-f)< {c_2}e^{-{t} {c_1}}\,.
\end{equation}
\end{theorem}
In the commutative case $f=1_A$\,, for some measurable subset $A$
of $\Omega$\,. The condition $\|(x-x_0)f\|\leq t$ implies that
$A\subseteq \{\omega\in \Omega:|x(\omega)-x_0(\omega)|\leq t\}$\,.
Then (\ref{eq:ldiv}) yields (\ref{eqlargedi}).

Our paper is organized as follows. In section 2, we explain the
notation and Kosaki's interpolation results, a fundamental tool in
our argument. In section 3, we discuss various inclusions of $BMO$
into $L_p$\,. Note that even the inclusion $BMO\subset L_2$ is not
immediate in the non-tracial setting. We first establish that for
$2\leq p< \infty\,, \ (BMO, L_p)$ forms an interpolation couple
and prove that $[BMO, L_p]_{\frac{p}{q}}=L_q\,, 2\leq p< q<
\infty$ (extending the results in \cite{Mu} to the non-tracial
setting). By general theory of interpolation, the continuous
inclusion $L_q\subset L_p$ yields a continuous inclusion
$BMO\subset L_p$\,. The inequality (\ref{JN9}) and the
modifications required to deduce the interval version of the
John-Nirenberg theorem are discussed in section 4, as well as the
proof of Theorem \ref{noncomldiv}\,. We are indebted to Tao Mei
and Narcisse Randrianantoanina for providing us with the preprints
\cite{Mei}, respectively \cite{Ra2}.

\section{Preliminaries}

We use standard notation in operator algebras. We refer to
\cite{KR} and \cite{Ta} for background on von Neumann algebras. In
the following, we will consider a $\sigma$-finite von Neumann
algebra $\N$ acting on a Hilbert space $H$\,, and a distinguished
normal faithful state $\phi$ on $\N$\,. We denote by
$\si_t=\si_t^\phi$ the one parameter modular automorphism group on
$\N$ associated with $\phi$\,. Haagerup's abstract $L_p$-spaces
are defined using the crossed-product
$\mathcal{R}=\N\rtimes_{\si_t}\rz$\,. We recall that
$\mathcal{R}\subseteq \mathcal{B}(L_2(\mathbb{R}, H))$ is the von
Neumann algebra generated by the operators $\pi(x)\,, x\in \N$
and, respectively, $\lambda(s)\,, s\in \mathbb{R}$\,, where
\[ \pi(x)(\xi)(t)=\si_{-t}(x)\xi(t) \quad\mbox{and}\quad
\la(t)(\xi)(s)=\xi(t-s)\,, \] for all $\xi\in
\mathcal{B}(L_2(\mathbb{R}, H))$ and all $t\in \mathbb{R}$\,. In
the following, we may and will identify $\N$ with $\pi(\N)$\,,
since $\pi$ is a normal faithful representation of $\N$ on
$L_2(\mathbb{R}, H))$\,. Furthermore, note that $\pi$ is invariant
under the dual action
\[ \theta_s(x) \lel W(s)xW(s)^*\,, \quad s\in \mathbb{R}\,, x\in \mathcal{R}\,, \]
and, moreover,
\[ \pi(\N)=\{x\in \mathcal{R}: \theta_s(x)=x\,, \,\,\text{for all}
\,\, s\in \mathbb{R}\}\,. \] Here the unitary operators $W(s)\,,
s\in \mathbb{R}$ are defined by the phase shift
\[ W(s)(\xi)(t) \lel e^{-ist}\xi(t)\,, \quad t\in \mathbb{R}\,.\]
As shown in \cite{PT}\,, the crossed product
$\N\rtimes_{\si_t}\rz$ is semifinite  and admits a unique normal
semifinite trace $\tau$ such that for all $s\in \mathbb{R}$\,,
 \[ \tau(\theta_s(x)) \lel e^{-s} \tau(x) \pl .\]
For $1\leq p\le \8$, $L_p(\N)$ is defined as the space of all
$\tau$-measurable operators $x$ affiliated with
$\mathcal{R}=\N\rtimes_{\si_t} \rz$ such that for all $s\in
\mathbb{R}$\,,
 \[ \theta_s(x)\lel e^{-{s}/{p}} x  \pl .\]
It follows from the definition that $L_\infty(\N)$ coincides with
$\N$\,. Furthermore, there is a canonical isomorphism between
$L_1(\N)$ and the predual $\N_*$ of $\N$\,. This requires some
explanation. Following \cite{PT}\,, every normal semifinite
faithful weight (n.s.f., for short) $\psi\in (\N_*)_+$ is given by
a density $h_\psi\in {L_1(\N)}_+$ satisfying
 \[ \tau(h_{\psi}x)=\int_{\mathbb{R}} \psi(\theta_s(x)) \ ds\,, \]
for all $x\in {\mathcal{R}}_+$\,. Using the polar decomposition of
an arbitrary element $\psi\in \N_*$\,, this correspondence between
$(\N_*)_+$ and ${L_1(\N)}_+$ extends to a bijection between $\N_*$
and $L_1(\N)$\,. Namely, if $\psi\in\N_*$\,, then
$\psi=u|\psi|$\,, where $u\in \N$ and $|\psi|$ is the modulus of
$\psi$\,. By construction, the corresponding $h_\psi\in {L_1(\N)}$
admits the polar decomposition
\[ h_\psi=u|h_\psi|=uh_{|\psi|}\,.
\] We may define a norm on $L_1(\N)$ by
\[ \|h_\psi\|_1=|\psi|(1)=\|\psi\|_{\N_*}\,, \quad \psi\in \N_*\,.
\]
In this way we obtain the isometry between $L_1(\N)$ and $\N_*$\,.
Furthermore, define a linear functional
$\text{tr}:L_1(\N)\rightarrow \mathbb{C}$\,, called {\em trace},
by
 \[ \text{tr}(h_{\psi}) \lel \psi(1)\,. \]
It is important to note that if ${\frac1p}+{\frac1q}=1$\,, then
for all $x\in L_p(\N)$ and $y\in L_q(\N)$\,, we have the tracial
property
 \[ tr(xy) \lel tr(yx)\,. \]
Let us now return to the special state $\phi$ fixed at the
beginning. Then, as explained above, there exists a density
$D_{\phi} \in L_1(\N)$ such that for all $x\in \N$\,,
 \[ \phi(x) \lel \text{tr}(D_{\phi}x)  \]
In the following we drop the subscript and reserve the letter $D$
exclusively for this density. Given $1\leq p< \infty$ and $x\in
L_p(\N)$\,, define
\[ \|x\|_p=(\text{tr}(|x|^p))^{1/p}\,, \quad
\text{respectively}\quad \|x\|_\infty=\|x\|_\N\,. \] If
${\frac1p}+{\frac1q}={\frac1r}$\,, then for all $x\in L_p(\N)$ and
$y\in L_q(\N)$\,, H\"{o}lder's inequality holds, i.e.,
\[ \|xy\|_r\leq \|x\|_p\|y\|_q\,. \]
As a consequence, given $1\leq p< \infty$\,, the mapping $(x,
y)\in {L_p(\N)\times L_{p'}(\N)}\mapsto \text{tr}(xy)\,,$ where
${\frac1p}+{\frac1{p'}}=1$\,, defines a duality bracket between
$L_p(\N)$ and $L_{p'}(\N)$\,, with respect to which we have the
isometry
\[ (L_p(\N))^*=L_{p'}(\N)\,. \]
In the sequel we will make repeated use of the following
well-known fact.
\begin{fact}\label{anfact} Let $1\le p,q,r\le \infty$ such that $\frac1p=\frac1q+\frac1r$ and $p<\infty$.
Let $x\in L_q(\N)$, and $a\in L_r(\N)$, $b\in L_r(\N)$ be positive
elements. Then the map
 \[ f(z)\lel a^{1-z}xb^{z} \in L_p(\N) \]
is continuous on $S=\{z\in \mathbb{C}\pl : \pl 0\le
\text{Re}(z)\le1\}$ and analytic in the interior $S_0$ of $S$.
\end{fact}

\noindent Consider now the 1-parameter automorphism group
 \[ \al_t(x)\lel D^{it}xD^{-it}\,, \quad t\in \mathbb{R}\,. \]
Note that $\al_t$ is strongly continuous and leaves $\N$ invariant
since
 \[ \theta_s(D^{it}xD^{-it})\lel e^{its}D^{it}xe^{-its} D^{its}
 \lel D^{it}xD^{-it} \pl .\]
From Fact \ref{anfact} we deduce that for every $x\in \N$\,, the
map $f(z)\lel D^{1-z}xD^z$ is analytic, taking values in
$L_1(\N)$\,. In particular, it follows that the map
$f_{x,y}(z)\lel \text{tr}(D^{1-z}xD^zy)$ is analytic and satisfies
 \[ f_{x,y}(it)\lel \phi(x\al_{t}(y)) \quad \mbox{and} \quad
 f_{x,y}(1+it)\lel \phi(\al_{t}(y)x) \pl ,\]
for all $t\in \rz$\,. Since
\[ \phi(\al_t(x))=\text{tr}(DD^{it}xD^{-it})=\phi(x)\,, \]
we deduce that $\si_t=\al_t\,,$ for all $t\in \mathbb{R}$ (see,
e.g., \cite[Section 2.12]{St}). \\
An element $x\in \N$ is called analytic if the map $t\in
\mathbb{R}\mapsto \si_t(x)\in \N$ extends to an analytic function
$z\in \mathbb{C}\mapsto \si_z(x)\in \N$. The family $\N_a$ of
analytic elements in $\N$ is a weak$^*$-dense $*$-subalgebra of
$\N$ (see \cite{PT})\,. Given $1\leq p< \infty$ and $0\leq
\theta\leq 1$\,, it can be shown (see, e.g., \cite{Ju}) that the
space ${\N_a}D^{1/p}$ is dense in $L_p(\N)$ and
\[ {D^{\frac{1-\theta}{p}}}{\N_a}{D^{\frac{\theta}{p}}}={\N_a}{D^{\frac1p}}\,. \]
We will denote by $L_p(\N)_a$ the family of elements with the
property that the map $t\mapsto \si_t(x),$ $t\in \mathbb{R}$
extends to an analytic function on $\cz$ with values in $L_p(\N)$.
Since $\N_aD^{\frac1p}\subset L_p(\N)_a$ is dense in $L_p(\N)\,,$
as mentioned above, it follows that $L_p(\N)_a$ is dense in
$L_p(\N)$.

Since interpolation is our main tool in the paper, we briefly
recall in the following some basic notions concerning the complex
method of interpolation due to Calder\'on. Our main reference for
interpolation theory is Bergh and L\"{o}fstr\"{o}m \cite{BL}. A
pair of Banach spaces $(X_0, X_1)$ is called a {\em compatible
couple} if they embed continuously in some topological vector
space $X$. This allows us to consider the spaces $X_0\cap X_1$ and
$X_0+X_1$ (by identifying them inside $X$). They are Banach spaces
when equipped, respectively, with the following norms
\begin{eqnarray}
\|x\|_{X_0\cap X_1}&=&\max\{\|x\|_{X_0}, \|x\|_{X_1}\}\,,\label{normonint}\\
\|x\|_{X_0+X_1}&=&\inf\{\|x_0\|_{X_0}+\|x_1\|_{X_1} :
x={x_0}+{x_1}\,, x_0\in X_0\,, x_1\in X_1\}\,.\label{normonsum}
\end{eqnarray}
Following the notation from \cite{BL},
let $\F$ be the family of all continuous and bounded functions
$f:S \rightarrow {X_0+X_1}$ satisfying  the following properties
\vspace*{0.2cm} \noindent
\begin{enumerate}
\item $f$ is analytic in $S_0$\,,\\[-0.2cm]
\item $f(it)\in X_0$ and $f(1+it)\in X_1$ \,for all \,$t\in {\mathbb{R}}$\,,\\[-0.2cm]
\item $f(it)\rightarrow 0$ and $f(1+it)\rightarrow 0$ \,as
\,$t\rightarrow \infty$\,.
\end{enumerate}
\vspace*{0.2cm} \noindent By the Phragmen-Lindel\"{o}f theorem,
$\F$ is a Banach space under the norm
\[ \|f\|_{\F} = \max \left\{\sup_{t \in {\mathbb{R}}}  \|f(it)\|_{X_0}, \;\sup_{t \in {\mathbb{R}}} \|f(1+it)\|_{X_1}\right\}\,. \]
For $0\leq \theta \leq 1$, set $[X_0,X_1]_{\theta}= \{\,x\in
{X_0+X_1} : x=f(\theta)\,, \:\mbox{for some} \:f\in {\F}\,\}$\,.
This is called the {\em complex interpolation space} (of exponent
$\theta$) between $X_0$ and $X_1$\,, and it is a Banach space
under the norm
\[ \|x\|_{\theta}= \inf \{\|f\|_{\F} : x=f(\theta), \: f\in {\F}\}\,. \]
The complex method is an exact interpolation functor; i.e., if
$T:{X_0}+{X_1}\rightarrow {{Y_0}+{Y_1}}$ is a linear operator
which is bounded both from $X_0$ to $Y_0$ (with norm $M_0$) and
from $X_1$ to $Y_1$ (with norm $M_1$), then $T$ is bounded from
$[X_0, X_1]_\theta$ to $[Y_0, Y_1]_\theta$\,, with norm $\leq
M_0^{1-\theta}M_1^\theta$\,.

\section{Interpolation results}
In this paper we will make crucial use of Kosaki's interpolation
results (see \cite{Ko})\,, which will enable us to extend the
results in \cite{Mu} to the non-tracial setting. For $1\leq p<
q\leq \infty$ and $0\leq \eta\leq 1$\, consider the mapping
$I_{q,p}^{\eta}:L_q(\N)\to L_p(\N)$ defined by
 \[ I_{q,p}^\eta(x)\lel D^{(1-\eta)(\frac1p-\frac1q)}xD^{\eta(\frac1p-\frac1q)} \pl .\]
Note that $I_{q,p}^{\eta}$ is injective. In \cite{Ko}, Kosaki
introduces the spaces
\begin{equation}\label{kos17}
 L_p^{\eta}(\N,\phi) \lel [I_{\infty,1}^{\eta}(\N),L_1(\N)]_{\frac
 1p} \pl ,
 \end{equation}
for $1\leq p\leq \infty$\,. We should point out that Kosaki's
original notation is slightly different, namely
$L_p^{\eta}(\N,\phi)$ is denoted in \cite{Ko} by
$i_p^\eta(L_p)$\,, where $i_p^\eta$ denotes the mapping $I_{p,
1}^\eta$\,. Also, in our terminology the space
$L_p^{\eta}(\N,\phi)$ is considered as a subspace of $L_1(\N)$
whereas Kosaki formally works in $\N_*$. Using the canonical
isometric isomorphism between $L_1(\N)$ and $\N_*$ explained in
the preliminaries, Kosaki's results may be reformulated in this
language as
\begin{equation}\label{kos24}
L_p^{\eta}(\N,\phi) \lel
 D^{\frac{1-\eta}{p'}}L_p(\N)D^{\frac{\eta}{p'}} \pl ,
 \end{equation}
 where $\frac{1}{p}+ \frac{1}{p'}=1\,.$
This is not only an equality on the level of sets, but it means
that if  $x\in L_p(\N)$ then
$D^{\frac{1-\eta}{p'}}xD^{\frac{\eta}{p'}}\in
[I_{\infty,1}^{\eta}(\N),L_1(\N)]_{\frac
 1p}$ \ and, furthermore,
 \[ \noo x\rrm_p \lel \noo
 D^{\frac{1-\eta}{p'}}xD^{\frac{\eta}{p'}}\rrm_{
 [I_{\infty,1}^{\eta}(\N),L_1(\N)]_{\frac
 1p}}. \]
Moreover, every element in the interpolation space
$[I_{\infty,1}^{\eta}(\N),L_1(\N)]_{\frac1p}$ comes form an
element $x\in L_p(\N)$. In our context, we will also consider the
interpolation couple $(A_0, A_1)$\,, where
$A_0=I_{q,p}^\eta(L_q(\N))$ and $A_1=L_p(\N)$\,, with $1\le p<
q\le \infty$\,. We will use the notation
 \[ \noo I_{q,p}^{\eta}(x)\rrm_{A_0} \lel \noo
 x\rrm_q \pl. \]
This is to emphasize that the ambient topological vector space for
the interpolation couple $(I_{q,p}^\eta(L_q(\N)),L_p(\N))$  is
$L_p(\N),$ and the inclusion map of $L_q(\N)$ into $L_p(\N)$ is
given by $I_{q,p}^\eta$\,. Using (\ref{kos17}), (\ref{kos24}) and
the reiteration theorem for complex interpolation (\cite[Theorem
4.6.1]{BL}), we deduce for $1\leq p< s< q\leq \infty$ and $0<
\theta< 1$ with $\frac{1}{s}=\frac{1-\theta}{q}+\frac{\theta}{p}$
\ that
\begin{eqnarray}
 I_{s,p}^{\eta}(L_s(\N))&\lel &
 [I_{q,p}^\eta(L_q(\N)),L_p(\N))]_{\theta}\,.\label{Kosinterp}
 \end{eqnarray}
We will now use these interpolation spaces in the context of
BMO-type norms. Indeed, it is helpful to consider the following
vector-valued $\ell_\8$ space defined for sequences $(x_n)_{n\geq
0}$ in $L_p(\N)$ as follows
 \[ \noo (x_n)\rrm_{L_p^c(\N;\ell_\8)} \lel \noo \sup_n (x_n^*x_n) \rrm_{\frac p2}^{\frac12}  \pl ,\]
where the right hand side is to be understood in the sense of the
following suggestive notation introduced in \cite{Ju}. Namely, if
$1\leq q,q' \leq \infty\,$ with ${\frac1q}+{\frac{1}{q'}}=1\,$
then, for a finite sequence of positive elements in $L_p(\N)$, set
\begin{eqnarray}
\|\sup_n {x_n}\|_q&=&\sup\left\{\,\sum_{n\geq 0}\text{tr}(x_ny_n)
: \; y_n\geq0, \left\|\sum_{n\geq 0} y_n\right\|_{q'}\leq
1\,\right\}\,.\label{defnormqsup}
\end{eqnarray}
In the commutative case this coincides with the usual definition
of the norm in the $l_\infty$- valued $L_q$-space. Here the actual
supremum does not necessarily make sense as in the commutative
case and thus it is just a \underline{notation}. Note that
$L_p^c(\N,\ell_\8)$ is only well-defined for $p\ge 2$.

\begin{rem}\label{rem5654}
{\rm Let $1\leq p< q, s< \infty$ be such that $\frac1p={\frac1q}+
{\frac1s}$. Then, for all $0\leq a, b\in L_s(\N),$
\begin{enumerate}
\item[$1)$] $\sup\limits_m \|\sup\limits_{n\leq m} x_n\|_p =
\|\sup\limits_n x_n\|_p$\\
\item[$2)$] $\|\sup\limits_n a^{\frac12}x_n a^{\frac12}\|_p\leq
\|a\|_s
\|\sup\limits_n x_n\|_q$\\
\item[$3)$] $\|\sup\limits_n x_n^{\frac12}b x_n^{\frac12}\|_p\leq
\|b\|_\infty \|\sup\limits_n x_n\|_p$\,.
\end{enumerate}
}
\end{rem}

\begin{proof}
Item 1) is an immediate consequence of the definition
(\ref{defnormqsup}). Item 2) follows from H\"{o}lder's inequality,
namely,
\begin{eqnarray*}
\|\sup\limits_n
a^{\frac12}x_na^{\frac12}\|_p&=&\sup\left\{\sum\limits_{n\geq 0}
\text{tr}(a^{\frac12}x_na^{\frac12}y_n): y_n\geq 0,
\left\|\sum\limits_{n\geq 0} y_n\right\|_{p'}\leq 1\right\}\\
&=& \sup\left\{\sum\limits_{n\geq 0}
\text{tr}(x_na^{\frac12}y_na^{\frac12}): y_n\geq 0,
\left\|\sum\limits_{n\geq 0} y_n\right\|_{p'}\leq 1\right\}\\
&\leq &\|a\|_s\|\sup\limits_nx_n\|_q\,.
\end{eqnarray*}

Item 3) follows from definition (\ref{defnormqsup}),
\begin{eqnarray*}
\|\sup\limits_n x_n^{\frac12}b
x_n^{\frac12}\|_p&=&\sup\left\{\sum\limits_{n\geq 0}
\text{tr}(x_n^{\frac12}bx_n^{\frac12}y_n): y_n\geq 0,
\left\|\sum\limits_{n\geq 0} y_n\right\|_{p'}\leq 1\right\}\\
&=& \sup\left\{\sum\limits_{n\geq 0}
\text{tr}(bx_n^{\frac12}y_nx_n^{\frac12}): y_n\geq
0, \left\|\sum\limits_{n\geq 0} y_n\right\|_{p'}\leq 1\right\}\\
&\leq &\|b\|_\infty \|\sup\limits_n x_n\|_p\,,
\end{eqnarray*}
wherein we have used the fact that
\[ \text{tr}(bx_n^{\frac12}y_nx_n^{\frac12})\leq
\|b\|_\infty \text{tr}(x_n^{\frac12}y_nx_n^{\frac12})=\|b\|_\infty
\text{tr}(x_ny_n)\,.
\]\\[-1.4cm]
\end{proof}

The following lemma, which can be proved in the same way as in the
tracial case (see \cite{Mu}), will enable us to use interpolation
for these norms. \lz

\begin{lemma}\label{no1} For $2\leq p\leq \infty$, let ${\frac1q}+
{\frac{2}{p}}=1$. Then
\[ \|(x_n)\|_{L_p^c(\N;\ell_\8)} =
\sup\left\{ \kla \summ_{n\geq 0}
 \noo x_nv_n\rrm_2^2\mer^{\frac12} : \noo \summ_{n\geq 0}
 v_nv_n^*\rrm_q \kl 1\right\} \pl . \]
 Consequently, $L_p^c(\N;\ell_\8)$ is a normed space.
\end{lemma}\lz

In the following, we will restrict our attention to the following
set-up. Let $\N$ be a $\sigma$-finite von Neumann algebra and
$\phi$ a distinguished normal, faithful state on $\N$. Let
$\M=\N\otimes \mathcal{B}(l_2)$ and let $\text{Tr}$ denote the
usual trace on $\mathcal{B}(l_2)$. Then $\phi\otimes \text{Tr}$ is
a n.s.f. weight on $\M$. For $1\leq p\leq \infty$\,, we consider
the Haagerup $L_p$-spaces $L_p(\M)$ associated with this n.s.f.
weight.

Let $V$ be the space of all infinite matrices $[x_{i, j}]_{1\leq
i, j\leq \infty}$ with entries in $L_p(\N)$, equipped with the
topology of pointwise convergence in $L_p(\N)$. More precisely, a
sequence of matrices $[x_{i, j}^{(m)}]$ converges to a matrix
$[x_{i, j}]$ if and only if
\begin{equation*}
L_p(\N)-\lim\limits_{m\rightarrow \infty} x_{i, j}^{(m)}= x_{i,
j}\,,
\end{equation*}
for all $1\leq i, j\leq \infty.$ Then $V$ is a topological vector
space. For $1\leq p< q\leq \infty$ and $0\leq \eta\leq 1$ define a
mapping
\[ \hat{I}_{q, p}^{\eta}={I}_{q,
p}^{\eta}\otimes {\text{Id}}_{\B(l_2)}\,. \] Note that
$(\hat{I}_{q, p}^{\eta}(L_q(\M)), L_p(\M))$ forms an interpolation
couple, since both embed continuously into $V$.

\begin{cor}\label{cor5678}
Let $1\leq p< s< q\leq \infty$\,, $0< \theta< 1$ such that
$\frac1s={\frac{1-\theta}{p}}+ {\frac{\theta}{q}}$ and $0\leq
\eta\leq 1$\,. Then, completely isometrically, \begin{eqnarray}
{[\hat{I}_{q, p}^{\eta}(L_q(\M)), L_p(\M)]}_\theta &=& \hat{I}_{s,
p}^{\eta}(L_s(\M))\,. \label{intlpm}
\end{eqnarray}
\end{cor}

\begin{proof}
By \cite{Ju1}, the following complete isometries hold
\[ L_p(\M)=S_p[L_p(\N)] \quad \text{respectively}, \quad
L_q(\M)=S_q[L_q(\N)]. \] Define $A_0={I}_{q, p}^{\eta}(L_q(\N))$
and  $A_1=L_p(\N)$. By the reiteration result (\ref{Kosinterp}) it
follows that $A_\theta=I_{s, p}(L_s(\N))$, completely
isometrically. Furthermore, by \cite[Corollary 1.4]{Pi2}, we
obtain  a complete isometry
\begin{eqnarray}
{[\,S_q[A_0], S_p[A_1]\,]}_\theta &=&
S_s[A_\theta].\label{rel8990}
\end{eqnarray}
Note that for $2\leq v< \infty$, we have
\begin{eqnarray*}
\| \, [x_{i, j} D^{\frac1p-\frac1q}] \, \|_{S_v[I_{v,
p}(L_v(\N))]}&=& \| \, [x_{i, j}] \, \|_{S_v[L_v(\N)]}.
\end{eqnarray*}
Together with (\ref{rel8990}), this yields the assertion.
\end{proof}

For $1\leq p< \infty$, recall the spaces $L_p(\N; l_2^c)$
(respectively, their row version $L_p(\N; l_2^r)$), defined in
\cite{JX} (see also \cite{PX} for the tracial case) as the
completion of the family of finite sequences $a=(a_k)_{k\geq 0}$
in $L_p(\N)$ under the norm \[  \|a\|_{L_p(\N;
l_2^c)}=\left\|(\sum\limits_{k\geq 0} |a_k|^2)^{\frac12}\right\|_p
\quad\mbox{respectively,}\quad \|a\|_{L_p(\N;
l_2^r)}=\left\|(\sum\limits_{k\geq 0}
|a_k^*|^2)^{\frac12}\right\|_p\,.
\]
In particular, for every positive integer $n$, we consider the
finite dimensional versions
\[  \|a\|_{L_p(\N;
(l_2^n)^c)}=\left\|(\sum\limits_{k=0}^n
|a_k|^2)^{\frac12}\right\|_p \quad\mbox{respectively,}\quad
\|a\|_{L_p(\N; (l_2^n)^r)}=\left\|(\sum\limits_{k=0}^n
|a_k^*|^2)^{\frac12}\right\|_p\,.
\]
From the above discussion, it is clear how to view the spaces
$(\hat{I}_{q, p}^{\eta}(L_q(\N;\ell_2^c)), L_p(\N;\ell_2^c)),$
where $2\leq p<q\leq \infty$ as a compatible couple and we obtain
the following:

\begin{cor}\label{co678}
Let $1\leq p< s< q\leq \infty$\,, $0< \theta< 1$ such that
$\frac1s={\frac{1-\theta}{p}}+ {\frac{\theta}{q}}$ and $0\leq
\eta\leq 1$\,. Then, completely isometrically,
\begin{eqnarray}
{[\hat{I}_{q, p}^{\eta}(L_q(\N;\ell_2^c)),
L_p(\N;\ell_2^c)]}_\theta &=&
 \hat{I}_{s, p}^{\eta}(L_s(\N;\ell_2^c)),\label{intlpml2c}\\
{[\hat{I}_{q, p}^{\eta}(L_q(\N;\ell_2^r)),
L_p(\N;\ell_2^r)]}_\theta &=& \hat{I}_{s, p}^{\eta}(
L_s(\N;\ell_2^r)),\label{intlpml2r}
\end{eqnarray}
\end{cor}

\begin{proof}
Note that $L_p(\N; l_2^c),$ respectively $L_p(\N; l_2^r)$ are
1-complemented subspaces of $L_p(\M)$ (see \cite{JX} for details).
Therefore the conclusion follows from Corollary \ref{cor5678}\,.
\end{proof}

\begin{rem}\label{ssweights}
{\rm
 Corollary \ref{cor5678} holds true, more generally, in the
setting of a von Neumann algebra $\M$ equipped with a n.f.
strictly semifinite weight $\psi$. Recall that a n.s.f. weight
$\psi$  is called \emph{strictly semifinite}, if there exists an
increasing sequence of projections $(p_n)_{n\geq 0}$ such that
$\si_t^{\psi}(p_n) \lel p_n$ and $\psi(p_n) \pl <\pl \8\,,$ for
all $n\geq 0$\,.
}
\end{rem}

We now establish additional properties of the spaces $L_p^c(\M;
l_\infty)$\,. Let $W$ be the space of infinite sequences of
matrices $[x_{i, j}^{(n)}]_{1\leq i, j, n\leq \infty}$ in
$L_p(\M)$\,. If $2\leq p< q\leq \infty$, then both $L_p^c(\M;
l_\infty)$ and $L_q^c(\M; l_\infty)$ embed continuously into the
topological vector space $W$. In fact, the map
${{\widehat{I}}}_{q, p}:L_q^c(\M; l_\infty)\rightarrow L_p^c(\M;
l_\infty)$ defined by
\[
{{\widehat{I}}}_{q, p}([x_{i, j}^{(n)}])=[x_{i, j}^{(n)}
D^{\frac1p-\frac1q}]
\] is a continuous inclusion from $L_q^c(\M; l_\infty)$ into $L_p^c(\M; l_\infty)$ and therefore we can interpolate between these
spaces.

The following result is probably well-known in interpolation
theory. We include a proof for the convenience of the reader.

\begin{lemma}\label{appint} Let $(A_0,A_1)$ be an interpolation couple with
intersection $\Delta$. Let $T_n:\Delta\to \Delta$ be a sequence of
linear operators such that
\begin{enumerate}
\item[1)] Each $T_n$ extends to a bounded linear operator from
$A_0$ to $A_0$ such that $C_1=\sup\limits_n \noo
 T_n\rrm<\infty$\,, and
\item[2)] Each $T_n$ extends to a bounded linear operator from
$A_1$ to $A_1$ such that $A_1-\lim\limits_{n\rightarrow \infty}
 T_n(x)\lel x\,,$ for all $x\in \Delta$\,.
\end{enumerate}
Then, for all $0< \theta< 1$ and all $x\in A_{\theta}=[A_0,
A_1]_\theta$\,, we have
\begin{equation}\label{eq45}
\lim\limits_{n\rightarrow \infty} T_n(x)\lel x\,.
\end{equation}
\end{lemma}
\begin{proof} By density \cite[Lemma 4.3.2]{BL} it suffices to prove
(\ref{eq45}) for elements $x\in A_\theta$ of the form
 \[ x\lel \summ_{k=0}^m g_k(\theta) y_k\,, \]
where $m$ is a non-negative integer, $g_k$ are analytic functions
which vanish at infinity and $y_k\in \Delta$\,. Let $\eps>0$ and
choose $n_0\geq 0$ such that for all $n\geq n_0$\,, we have
 \[ ((1+C_1)\summ_{k=0}^m \sup_t|g_k(it)| \noo y_k\rrm_{A_0})^{1-\theta} (\summ_{k=0}^m \sup_t |g_k(1+it)| \noo T_n(y_k)-y_k\rrm_{A_1})^{\theta} <\eps \pl .\]
Since the complex interpolation method is an exact functor of
exponent $\theta$\,, we deduce that
\begin{eqnarray*}
 \noo x-T_n(x)\rrm_\theta
 &\le &\sup\limits_{t\in \mathbb{R}} \noo \summ_{k=1}^m g_k(it) (y_k-T_n(y_k))\rrm_{A_0}^{1-\theta}
 \sup\limits_{t\in \mathbb{R}} \noo \summ_{k=1}^m g_k(1+it) (y_k-T_n(y_k))\rrm_{A_1}^{\theta} \\
 &\le &((1+C_1)\summ_{k=0}^m \sup_t |g_k(it)| \noo
 y_k\rrm_{A_0})^{1-\theta} (\summ_{k=0}^m \sup_t |g_k(1+it)|\noo
 y_k-T_n(y_k)\rrm_{A_1})^{\theta}\\
&< &\eps\,.
 \end{eqnarray*}
This proves the assertion.
\end{proof}

\begin{rem}\label{rem45}
{\rm Let $2\leq p< q\leq \infty$ and $e\in \M$ be a projection
such that $(\phi\otimes \text{Tr})(e)< \infty.$ Then $e{\M}e$ is a
finite von Neumann algebra, and a similar argument as in the
tracial case (see \cite{Mu}) shows that the following contractive
inclusion holds
\begin{eqnarray}
{{\widehat{I}}}_{q, p}(L_q^c(e{\M}e; l_\infty))&\subseteq &
L_p^c(e{\M}e; l_\infty)\,.\label{rel7867}
\end{eqnarray}
}
\end{rem}

\begin{prop}\label{inte} If $2\le p< s< q\leq \8$ and $0< \theta< 1$ such that
$\frac1s=\frac{1-\theta}{q}+\frac{\theta}{p}$\,, then,
isometrically,
\begin{eqnarray}\label{eq99}
{[{{\widehat{I}}}_{q, p}(L_q^c(\M;\ell_\8)),
L_p^c(\M;\ell_\8)]}_\theta &=&
 {{\widehat{I}}}_{s, p}(L_s^c(\M;\ell_\8)).
 \end{eqnarray}
\end{prop}\lz

\begin{proof}  We first show that we have a contractive
inclusion
\begin{eqnarray}
{[{{\widehat{I}}}_{q, p}(L_q^c(\M;\ell_\8)),
L_p^c(\M;\ell_\8)]}_\theta &\subseteq &
 {{\widehat{I}}}_{s, p}(L_s^c(\M;\ell_\8)).\label{rel9090}
 \end{eqnarray}
Note that for each positive integer $n$\,, $L_p^c(\M; l_\infty^n)$
is isomorphic (as a vector space) to the $n$-fold tensor product
${L_p(\M)}\otimes \ldots \otimes L_p(\M)$ ($n$ times). Therefore
$({{\widehat{I}}}_{q, p}^{(n)}(L_q^c(\M; l_\infty^n)), L_p^c(\M;
l_\infty^n))$ is a compatible couple, where
\[ {{\widehat{I}}}_{q, p}^{(n)}= \hat{I}_{q, p}\otimes
{\text{Id}(l_\infty^n)}. \]By Remark \ref{rem5654} (1), it
suffices to show that for all $n\geq 1$ we have, contractively,
\begin{eqnarray}
{[{{\widehat{I}}}_{q, p}^{(n)}(L_q^c(\M;\ell_\8^n)),
L_p^c(\M;\ell_\8^n)]}_\theta \subseteq
 {{\widehat{I}}}_{s, p}^{(n)}(L_s^c(\M;\ell_\8^n)).\label{rel4534}
\end{eqnarray}
In the following we fix a positive integer $n$ and consider only
sequences of length $n$. We claim that it suffices to prove that
for all $m\geq 1$ the following contractive inclusion holds
\begin{eqnarray}
{[{{\widehat{I}}}_{q, p}^{(n)}(L_q^c(\N\otimes M_m; l_\infty^n)),
L_p^c(\N\otimes M_m; l_\infty^n)]}_\theta &\subseteq &
{{\widehat{I}}}_{q, p}^{(n)}(L_s^c(\N\otimes M_m; l_\infty^n))\,,
\label{rel8899}
\end{eqnarray}
where $M_m$ is the algebra of $m\times m$ complex matrices. Note
that $\N\otimes M_m=e_m{\M} e_m$, where $e_m=1\otimes p_m$ and
$p_m: \mathcal{B}(l_2)\rightarrow M_m$ is the canonical orthogonal
projection. Set \[ A_0={{\widehat{I}}}_{q,
p}^{(n)}(L_q^c(\N\otimes M_m; l_\infty^n)) \quad\mbox{and} \quad
A_1=L_p^c(\N\otimes M_m; l_\infty^n)\,.
\] By Remark \ref{rem45}, their intersection and sum are,
respectively,
\[ \Delta_0 = A_0\cap A_1= {{\widehat{I}}}_{q, p}^{(n)}(L_q^c(\N\otimes M_m; l_\infty^n)), \quad \Si_0    =A_0+A_1= L_p^c(\N\otimes M_m ;\ell_\8^n)\,. \]
Consider the map $S_m: L_q^c(\M; l_\infty^n)\rightarrow
L_q^c(e_m{\M}e_m; l_\infty^n)$ defined by
\[ S_m((x_k)_{k=1}^n)=(e_m{x_k}e_m)_{k=1}^n. \]
By (2) in Remark \ref{rem5654}, it follows that
\[ \sup\limits_m\|S_m\|_{A_0\rightarrow A_0}\leq 1\,. \]
Furthermore, we have $A_1-\lim\limits_{m\rightarrow \infty}
S_m(x)=x$, for all $x\in \Delta_0$ (see \cite{Ju}). An application
of Lemma \ref{appint} shows that, for all $\ x\in
[A_0,A_1]_\theta$\,,
\[ \lim\limits_{m\rightarrow \infty} S_m(x)=x\,. \]
Therefore the claim is justified. We now prove (\ref{rel8899}).
Let $2\leq v< \infty$ such that $\frac1v=\frac1p-\frac1q\,.$
Furthermore, define $f(v)$ by the relation $\frac1{f(v)}+
\frac1v=\frac12$ and note that
$\frac{1}{f(q)}-\frac{1}{f(p)}=\frac1v$\,. Let
\[ B_0={{\widehat{I}}}_{f(p), f(q)}^{(n)}(L_{f(p)}(\N\otimes M_m; (l_2^n)^r))\,,
\quad B_1={{\widehat{I}}}_{f(q), f(q)}^{(n)}(L_{f(q)}(\N\otimes
M_m; (l_2^n)^r))\,.
\] Then, their intersection and sum are, respectively, \[
\Delta_1={{\widehat{I}}}_{f(p), f(q)}^{(n)}(L_{f(p)}(\N\otimes
M_m; (l_2^n)^r))\,, \quad \Si_1 \lel
 {{\widehat{I}}}_{f(q), f(q)}^{(n)}(L_{f(q)}(\N\otimes M_m; (l_2^n)^r ))\,. \]
Furthermore, given finite sequences of $m\times m$ matrices
$([x_{i,j}^{(k)} D^{\frac1v}]_{1\leq i, j\leq m})_{k=1}^n$ in
$\Delta_0$ and, respectively, $([D^{\frac1v}y_{i,j}^{(k)}]_{1\leq
i, j\leq m})_{k=1}^n$ in $\Delta_1$, define
 \[ T\left(\left([x_{i, j}^{(k)}D^{\frac1v}]\right)_k\times  \left([D^{\frac1v}y_{i,j}^{(k)}]\right)_k\right) \lel
 \left([x_{i, j}^{(k)}D^{\frac1v} y_{i, j}^{(k)}]\right)_{k} \in \ell_2^n(L_2(\N\otimes M_m)) \pl .\]
This map is well-defined, since for each $1\leq k\leq n$, $[x_{i,
j}^{(k)}]_{1\leq i, j\leq m} \in L_q(\N\otimes M_m)$ and $[y_{i,
j}^{(k)}]_{1\leq i, j\leq m}\in L_{f(p)}(\N\otimes M_m)\,,$ and,
by construction,  $\frac1q +\frac1v + \frac1{f(q)}=\frac12\,.$
Moreover, we deduce from Lemma \ref{no1} that
 \[ \noo T: \Delta_0 \times \Delta_1 \to \ell_2^n(L_2(\N\otimes M_m))\rrm \kl 1 \pl .\]
Using multilinear complex interpolation (see \cite[Theorem
4.4.1]{BL}), the map $T$ extends to a contraction
\[ T:{[A_0, A_1]_\theta}\times {[B_0, B_1]_\theta}\to \ell_2^n(L_2(\N\otimes
M_m))\,. \] By Corollary \ref{co678} applied to the algebra
$\N\otimes M_m$\,, we deduce that, isometrically,
\[ {[B_0, B_1]_\theta}= {{{\widehat{I}}}_{f(s), f(q)}^{(n)}(L_{f(s)}(\N\otimes M_m;
(l_2^n)^r
))}\,.
\] Therefore the extended map
\[ T: [A_0, A_1]_\theta \times
  {{\widehat{I}}}_{f(s), f(q)}^{(n)}(L_{f(s)}(\N\otimes M_m; (l_2^n)^r ))\to \ell_2^n(L_2(\N\otimes
  M_m)) \]
is a contraction. A further application of Lemma \ref{no1} yields
the estimate
 \begin{equation*} \label{TC}
  \noo T:[{{\widehat{I}}}_{q, p}^{(n)}(L_q^c(\N\otimes M_m; \ell_\8^n)), L_p^c(\N\otimes M_m; \ell_\8^n)]_\theta  \to {{\widehat{I}}}_{s, p}^{(n)}(L_s^c(\N\otimes M_m;\ell_\8^n)) \rrm \kl
  1 \pl .
 \end{equation*}
Thus (\ref{rel8899}) is justified. The proof of the reverse
inclusion in (\ref{eq99}) reduces, just as in the tracial case, to
the interpolation result (\ref{intlpm}), via the factorization
property of the spaces $L_p^c(\M; l_\infty),$ proved in \cite{Ju}.
\end{proof}

Let us now turn our attention to  martingales invariant under the
automorphism group of a given state. Indeed, we will assume that
$(\N_n)_{n\geq 0}$ is an increasing sequence of von Neumann
subalgebras of $\N$ such that $\si_t(\N_n) \subset \N_n$ for all
$t\in \mathbb{R}$ and all $n\geq 0$\,. Moreover, let $(p_n)_{n\geq
0}$ be the increasing sequence of projections given by the units
in $\N_n$\,. We will further assume that $\si_t(p_n)=p_n$ for all
$t\in \mathbb{R}$ and all $n\geq 0$\,. For fixed $\nen$\,,
according to \cite{Ta1}\,, there exists a conditional expectation
$\tilde{\E}_n:p_n\N p_n\to \N_n$ such that
 \begin{equation}
 \label{inv} \tilde{\E_n}(\si_t(x)) \lel \si_t(\tilde{\E_n}(x))
 \end{equation}
for all $x\in p_n \N p_n$ and all $t\in \mathbb{R}$\,. Let us
define $\E_n(x)\lel \tilde{\E_n}(p_nxp_n)$ and note that $\E_n$
still satisfies the invariance property (\ref{inv}). Then the
martingale differences are defined by
 \[ d_n(x) \lel  \E_{n}(x)-\E_{n-1}(x) \pl ,\]
where $\E_{-1}(x)=0$. Moreover, for $1\leq p< \infty$\,, the
conditional expectation and the martingale difference operators
extend in  a natural way to $L_p(\N)$ (see \cite{JX}),  such that
for all $0< \theta< 1$ the relation
 \begin{equation} \label{exp}
 \E_n(D^{\frac{1-\theta}{p}}xD^{\frac{\theta}{p}})
 \lel
  D^{\frac{1-\theta}{p}}\E_n(x)D^{\frac{\theta}{p}}
  \end{equation}
holds for all $x\in \N$. Furthermore, if
$\frac1p={\frac1q}+{\frac1s}+{\frac1t}$\,, then for all $a\in
L_s(\N_n)$, $b\in L_t(\N_n)$ and $x\in L_q(\N)$\,, we have
 \begin{equation} \label{module}  \E_n(axb) \lel a\E_n(x)b\,.
  \end{equation}
Recall also that if $1\leq p\leq \infty$ and $x_\infty\in
L_p(\N)$\,, then the sequence $(x_n)_{n\geq 0}$ defined by
$x_n=\E_n(x_\infty)$ is a bounded $L_p(\N)$-martingale  which
converges to $x_\infty$ in $L_p(\N)$ (respectively, in the
w$^*$-topology if $p=\infty$)\,. Conversely, for $1< p< \infty$\,,
we deduce from the uniform convexity of $L_p(\N)$ that every
bounded $L_p(\N)$-martingale converges to some element
$x_\infty\in L_p(\N)$ and hence it is of the above form.
Therefore, we will often identify a martingale with its limit,
whenever this exists.

For convenience, we will also assume that the predual $\N_*$ is
separable. Our results will hold without this assumption and we
refer to the Appendix in \cite{Gous} for techniques deducing the
general case from the $\si$-finite one. However, assuming that
$\N_*$ is separable, we will find a separable $\si$-weakly dense
subalgebra $\A_0$ of $\N_a$. Following \cite{Ju}, we can then
consider the countably generated Hilbert $C^*$-module $F_n$
generated by $\N_n$ and $\A_0$\,, and find a right module
isomorphism $u_n:F_n\to C_\8(\N_n)$ such that for all $x, y\in
F_n$\,, we have
\begin{eqnarray}
\E_n(y^*x)&=&u_n(y)^*u_n(x)\,.\label{rel89907}
\end{eqnarray}
Moreover, for $1\leq p<\8$ the map $u_n^p: \N_{n,a}D^{\frac1p}\to
C_\8(\N_{n,a})D^{\frac1p}$ defined by
\[ u_n^p(xD^{\frac1p}) \lel u_n(x)D^{\frac1p}\,, \quad x\in {\N_{n,a}} \]
is well-defined and isometric if we equip $\N_{n,a}D^{\frac1p}$
with the norm
 \[ \noo xD^{\frac1p} \rrm_{L_p(\N,\E_n)}  \lel \noo
 D^{\frac1p}\E_n(x^*x) D^{\frac1p}
 \rrm_{\frac{p}{2}}^{\frac12} \pl .\]
Indeed, in \cite{Ju} this was considered only in the faithful
case. However, by working  with the algebra $\N_n+\cz(1-p_n)$ and
the corresponding $\tilde{u}_n^p\,,$ and then defining
$u_n^p(xD^{\frac1p})=\tilde{u}_n^p(xD^{\frac1p})p_n$, we can
extend the construction to the non-faithful case (see Section 8 in \cite{JX} for details).\\
For $2\leq p\leq \infty$\,, we will consider the space
$L_p^cMO(\N)$ of martingale difference sequences $x\cong
(d_k)_{k\geq 0}$ satisfying
\begin{equation}\label{lpcmonorm}
 \noo x\rrm_{L_p^cMO(\N)} \lel \sup_{m}   \noo
 \sup_{n \le m} \E_n\left(\summ_{k=n}^m d_k^*d_k\right)
 \rrm_{\frac p2}^{\frac12}< \infty \pl .
\end{equation}
Note that if $p=\infty$\,, then
${L_\infty^cMO(\N)}=BMO_c(\N)$\,.\\
Recall also the following norms
\begin{equation}\label{eq7879}
\|x\|_{H_p^c(\N)}=\|dx\|_{L_p(\N; l_2^c)}\,.
\end{equation}
The next results concerning the spaces $L_p^cMO(\N)$ will be very
useful in the sequel. The following noncommutative version of
Stein's inequality is essentially contained in \cite{JX}.

\begin{lemma}\label{Stbmo}
Let $1< p< \infty$\,. Define a map $Q$ on all finite sequences
$x=(x_n)_{n\geq 0}$ in $L_p(\N)$ by $Q(x)=(\E_n(x_n))_{n\geq
0}$\,. Then
\[ \|Q(x)\|_{L_p(\N; l_2^c)}\leq \gamma_p\|x\|_{L_p(\N;
l_2^c)} \quad\text{respectively}\,, \quad \|Q(x)\|_{L_p(\N;
l_2^r)}\leq \gamma_p\|x\|_{L_p(\N; l_2^r)}\,. \] Thus $Q$ extends
to a bounded projection on $L_p(\N; l_2^c)$ and $L_p(\N; l_2^r)$,
respectively.
\end{lemma} \lz

Furthermore, as an application of Stein's inequality, Doob's
inequality and the duality between $L_p^cMO(\N)$ and
$H_{p'}^c(\N)$\,, where $\frac1p+\frac1{p'}=1$\,, the following
was proved in \cite{Ju}(see also Section 8 in \cite{JX}):

\begin{theorem} \label{1step}
Let $1< p < \8$, then, with equivalent norms,
 \[ H_p^c(\N)= L_p^cMO(\N)\,. \]
\end{theorem}\lz

Let $2\leq p< q\leq \infty$. As already observed in \cite{JX}, the
map $I_{q,p}:L_q^cMO(\N)\to L_p^cMO(\N)$ is a contraction.
However, for $0< \eta< 1$, the inclusion map ${I_{q,p}^\eta}:
L_q^cMO(\N)\to L_p^cMO(\N)$ is not necessarily continuous anymore,
but well-defined on finite sequences. Since we would like to
consider interpolation between these spaces, let
$\mathcal{S}=\prod\limits_{n\in \mathbb{N}}\prod\limits_{n\le m}
L_p(\N_m)$, endowed with the product topology. For $p< v< q$ with
$\frac1p=\frac1q+\frac1v$, we may consider the embedding
$\tilde{I}_{q,p}^{\eta}:L_q^cMO(\N)\to \mathcal{S}$  defined by
\begin{equation}\label{e.3.333}
\tilde{I}_{q,p}^{\eta}((x_m)_m) \lel
 (D^{\frac{1-\eta}{v}}(x_m-x_{n-1})D^{\frac{\eta}{v}})_{m, n}
\end{equation}
where $(x_m)_m$ is a finite martingale. Note that by theorem
\ref{1step}, it follows that finite martingales are dense in
$L_q^cMO(\N)$ for $q< \infty$\,, since this is true for the space
$H_q^c(\N)$\,.

The following proposition is a non-tracial modification of the
result in \cite{Mu}.
\begin{prop}\label{intlqcmo}  Let $2\le p < s< q\le \8$ and $0\le \theta \le 1$ such that $
\frac1s=\frac{1-\theta}{q}+ {\frac{\theta}{p}}$, then, with
equivalent norms
\begin{eqnarray}
[\tilde{I}_{q,p}(L_q^cMO(\N)), L_p^cMO(\N)]_\theta
&=&\tilde{I}_{s,p}(L_s^cMO(\N))\,.\label{lpcmointerp}
\end{eqnarray}
\end{prop}\lz

\begin{proof}
For $2\le v\le \8$, we define a map $\Phi_v: L_v^cMO(\N)\to
L_v^c(\M;\ell_\8)$ by
 \[ \Phi_v(x)=( u_n^v(x-{\mathcal{E}}_{n-1}(x)))_{n\geq 0}  \pl .\]
We claim that $\Phi_v$ is an isometry. Indeed, an application of
Lemma \ref{no1}, together with (\ref{rel89907}) shows that
 \begin{eqnarray*}
 \noo \Phi_v(x)\rrm_{L_v^c(\M;\ell_\8)}
  &=& \left\|\sup\limits_n u_n^v(x-{\mathcal{E}}_{n-1}(x))^*u_n^v(x-{\mathcal{E}}_{n-1}(x))\right\|_{\frac{v}{2}}^{\frac12} \\
  &=& \left\|\sup\limits_n {\mathcal{E}}_n((x-{\mathcal{E}}_{n-1}(x))^*(x-{\mathcal{E}}_{n-1}(x))) \right\|_{\frac{v}{2}}^{\frac12} \\
  &=& \noo x\rrm_{L_v^cMO(\N)} \pl .
  \end{eqnarray*}
Moreover,  we observe that by definition of the maps $u_n^v$ we
have
 \[  \Phi_p(xD^{\frac1p-\frac1s}) \lel
 \Phi_s(x)D^{\frac1p-\frac1s} \pl .\]
This implies that
 \[ \Phi_p \tilde{I}_{s, p} \lel {{\widehat{I}}}_{s, p}\Phi_s\,. \]
Hence, the family $\Phi_v$ is compatible with the interpolation
couples
 \[ (\tilde{I}_{q,p}(L_q^cMO(\N)), L_p^cMO(\N))\quad \mbox{and} \quad
  ({{\widehat{I}}}_{q, p}(L_q^c(\M,\ell_\8)), L_p^c(\M,\ell_\8))\pl .\]
We deduce that
 \[ \Phi_p: [\tilde{I}_{q, p}(L_q^cMO(\N)), L_p^cMO(\N)]_\theta \to
 [{{\widehat{I}}}_{q, p}(L_q^c(\M;\ell_\8)), L_p^c(\M;\ell_\8)]_\theta \]
is a contraction. An application of Proposition \ref{inte} yields
the contraction
 \[ \widehat{I}_{s, p}\Phi_p: [\tilde{I}_{q, p}(L_q^cMO(\N)), L_p^cMO(\N)]_\theta \to
 \widehat{I}_{s, p}(L_s^c(\M;\ell_\8)) \pl. \]
By considering elements in the intersection
$\Delta=(L_q^cMO(\N))D^{\frac1p-\frac1q}\subset L_p^cMO(\N)$, we
deduce that the image of ${{\widehat{I}}}_{s, p}\Phi_p$ is
contained in ${{\widehat{I}}}_{s, p}\Phi_s(L_s^cMO(\N))$. Since
$\Phi_s$ is isometric, we obtain a contractive inclusion
\begin{eqnarray}
 j_s: [\tilde{I}_{q, p}(L_q^cMO(\N)), L_p^cMO(\N)]_\theta \to \tilde{I}_{s,
 p}(L_s^cMO(\N))\,,\label{eq55558}
 \end{eqnarray}
densely defined  by the formula
 \[ j_s(xD^{\frac1p-\frac1q}) \lel
 xD^{\frac1p-\frac1s} \pl .\] A similar argument as in the proof of Proposition 3.4 in \cite{Mu} yields the converse inclusion. Namely, since the unit ball of the interpolation space $[\tilde{I}_{q,
p}(L_q^cMO(\N)), L_p^cMO(\N)]_\theta$ is dense in the unit ball of
the space $[\tilde{I}_{q, p}(L_q^cMO(\N)), L_p^cMO(\N)]^\theta$ with
respect to the sum topology, we also obtain the contractive
inclusion $[\tilde{I}_{q, p}(L_q^cMO(\N)), L_p^cMO(\N)]^\theta
\subseteq \tilde{I}_{s, p}(L_s^cMO(\N))$\,. Therefore, to show the
converse inclusion in (\ref{eq55558}) it follows by duality that it
suffices to prove
\begin{eqnarray}
[H_{q'}^c(\N), \tilde{I}_{p', q'}(H_{p'}^c(\N))]_\theta \subseteq
\tilde{I}_{s', q'}(H_{s'}^c(\N))\,,\label{eq333356}
\end{eqnarray}
where $p', q', s'$ are respectively the conjugates of $p, q, s$\,.
Note that since $\frac1s=\frac{1-\theta}{q}+ {\frac{\theta}{p}}$\,,
we deduce that $ 1\leq q'< s'< p'\leq 2$ and
$\frac1{s'}=\frac{1-\theta}{q'}+ {\frac{\theta}{p'}}$\,. Therefore,
by Corollary \ref{co678} we obtain
\begin{eqnarray}
[H_{q'}^c(\N), \tilde{I}_{p', q'}(H_{p'}^c(\N))]_\theta\subseteq
{[L_{q'}(\N; l_2^c), \hat{I}_{p', q'}(L_{p'}(\N; l_2^c))]}_\theta
= \hat{I}_{s', q'}(L_{s'}(\N; l_2^c))\,.\label{eq3334563}
\end{eqnarray}
By Stein's inequality (see Lemma \ref{Stbmo}) we have for all
$x=(x_m)_m\in H_{q'}^c(\N)\subseteq \tilde{I}_{s',
q'}(H_{s'}^c(\N))$\,,
\[ \left\|\left(\sum d_m(x_m)^*d_m(x_m)\right)^\frac12\right\|_{s'}\leq
\gamma_{s'}\left\|\left(\sum x_m^*
x_m\right)^\frac12\right\|_{s'}\,,
\] where the constant $\gamma_{s'}$ depends only on $s'$\,. Together
with (\ref{eq3334563}), this implies (\ref{eq333356}) and the
assertion is proved.
\end{proof}

As an application, we prove the result announced in \cite{JX}.
\lz

\begin{cor}\label{corapplic} Let $4< p < \8$, then  for every $x\in
L_p(\N)$\,, we have
 \[ \noo x\rrm_{L_p^cMO(\N)} \kl c_4^{{4}/{p}} \pl  \noo x\rrm_p \pl .\]
Here $c_4$ is an absolute constant.
\end{cor} \lz

\begin{proof} Fix $0<\theta < 1$ such that
$\frac1p=\frac{1-\theta}{\8}+\frac{\theta}{4}$\,.  Let $x\in
L_p(\N)$ with $\|x\|_p\leq 1\,.$ Consider the polar decomposition
$x=u|x|$ with $u\in \N$ and $|x|=(x^*x)^{\frac12}\in L_p(\N)$\,.
It follows that $\noo |x|\rrm_p =\noo x\rrm_p\le 1$. Given
$\delta>0$ and $C\ge 1$ define a function
 \[ f(z) \lel  C^{-z} \exp(\delta(\theta-z)^2) \pl u |x|^{\frac{zp}{4}}D^{\frac{(1-z)}{4}} \pl ,\]
for all $z$ in the unit strip. According to Fact \ref{anfact} in
the preliminaries, it follows that $f$ is analytic and $f(z)\in
L_4(\N)$\,. Furthermore, for all $t\in \rz\,,$ note that
$|x|^{pit}D^{-it}\in \N\,.$ Hence
\begin{equation}\label{eq3454}
 \noo f(it) D^{-{\frac14}}\rrm_\8 \kl \exp(\delta \theta^2)\,,
 \end{equation}
and, respectively,
\begin{equation}\label{eq3455}
 \noo f(1+it) \rrm_4  \kl  C^{-1} \pl \exp(\delta \theta^2) \pl .
 \end{equation}
According to Theorem \ref{1step} (see \cite{JX} for details),
there exist absolute constants $\beta_4\,, \gamma_4> 0$ such that
for all $y\in L_4(\N)\,,$
 \[ \noo y\rrm_{L_4^cMO(\N)} \kl 2\gamma_4 \pl \noo y\rrm_{H_4^c(\N)} \kl
 2\gamma_4 \beta_4 \noo y\rrm_4 \pl .\]
Trivially, we have for $y\in \N$ that
 \[ \noo y\rrm_{BMO_c(\N)} \kl 2\noo y\rrm_\8 \pl .\]
Choose $C=2\gamma_4 \beta_4:=c_4$\,. It follows that $f(\theta)\in
[\tilde{I}_{\infty, 4}(BMO_c(\N)), L_4^cMO(\N)]_{\theta}\,,$ with
\[ \|f(\theta)\|_\theta\leq
(\exp(\delta\theta^2))^{1-\theta}(\exp(\delta\theta^2))^{\theta}=\exp(\delta)\,.
\]
According  to Proposition \ref{intlqcmo}, we deduce that
 \[ \noo \tilde{I}_{p, 4}(f(\theta))\rrm_{L_p^cMO(\N)} \kl \exp(\delta) \pl .\]
Since $\delta>0$ was arbitrarily chosen and
$f(\theta)={c_4}^{-\frac{4}{p}}xD^{\frac14-\frac1p}$, the
assertion follows.
\end{proof}

\begin{lemma}\label{densarg}  Let $2\leq p<q\le \infty$\,, $0<\theta<1$ and $0\leq \eta\leq 1$. Then
$\bigcup\limits_{m \geq 0} \tilde{I}_{q,p}((\N_{m,a})D^{\frac1q})$
is norm dense in the interpolation space $\
[\tilde{I}_{q,p}^{\eta}(L_q^cMO(\N)), L_p^cMO(\N)]_{\theta}$.
\end{lemma}

\begin{proof} First we note that the sequence $(\E_m)_{m\geq 0}$ given by the conditional
expectations satisfies the assumptions of Lemma \ref{appint}. Thus
the space of finite martingales is norm dense in the interpolation
space, since this is true for the space $L_p^cMO(\N)$\,, as
observed before. Given an integer $m\geq 0$\,, we note that the
inclusions
 \[ L_s(\N_m; l_2^c) \subset L_s^cMO(\N_m) \subset L_s^cMO(\N) \]
are continuous, since
 \[ \noo d_k(x)\rrm_{L_s(\N_m)}\kl 2 \min\{\noo
 x\rrm_{L_s^cMO(\N_m)},\noo x\rrm_{L_s(\N_m)}\} \pl \]
for all $0\le k\le m$. The intersection space is
\[ \tilde{I}_{q,p}^{\eta}(L_q^cMO(\N_m))\cap L_p^cMO(\N_m)\lel
\tilde{I}_{q,p}^{\eta}(L_q^cMO(\N_m))\,. \] Since
$\N_mD^{\frac1q}$ is norm dense in $L_q(\N_m)$, we deduce that
$\tilde{I}_{q, p}^{\eta}(\N_mD^{\frac1q})$ is also dense in the
interpolation space. In order to obtain analytic elements we
recall the approximation operators $R_k:L_p(\N)\to L_p(\N)$
defined by
 \[ R_k(x) \lel k^{\frac12} \pi^{-\frac12} \int_{\rz} \exp(-kt^2) D^{it}xD^{-it} \pl dt \pl .\]
Since $\si_t(x)=D^{it}xD^{-it}$ is strongly continuous on $\N$, we
deduce by approximation that $\si_t$ is norm continuous on
$L_p(\N)$ for all $2\leq p<\infty$ (see \cite{Ju}, \cite{JS}). It
then follows that
 \[ L_p-\lim\limits_{k\rightarrow \infty} R_k(x)\lel x \]
for every $x\in L_p(\N)$.  Therefore, the same is true for the
space $L_p^cMO(\N_m)$. Clearly, the map $\si_t:L_p^cMO(\N)\to
L_p^cMO(\N)$ is a contraction. It follows that $R_k$ is a
contraction on $L_p^cMO(\N)\,,$ for all $2\le p\le \infty$. In
view of Lemma \ref{appint}, we deduce that elements of the form
$\tilde{I}_{q, p}^{\eta}(R_k(x)D^{\frac1q})$ are norm dense in the
interpolation space $[\tilde{I}_{q,p}^{\eta}(L_q^cMO(\N_m)),
L_p^cMO(\N_m)]_{\theta}$\,, for all $x\in \N_m$\,. Using the fact
that the image $R_k(\N_m)$ consists of analytic elements and
 \[ R_k(xD^{\frac1q})\lel R_k(x)D^{\frac1q} \pl ,\]
we deduce the assertion.
\end{proof}

In order to establish the inclusion $I_{\infty,
p}^1(BMO(\N))=I_{\infty, p}^1({BMO}_c\cap {BMO}_r)\subset
L_p(\N)$\,,  we have to consider and analyze the space $I_{\infty,
p}^1({BMO}_r(\N))$\,. In contrast with the inclusion $I_{\infty,
p}^1({BMO}_c(\N))$ $\subset L_p^c{MO}(\N)$\,, which follows easily
by the definition, for the row version we are multiplying by
$D^{\frac1p}$ on the wrong side. Alternatively, we may consider
the space $I_{\infty, p}^0({BMO}_c(\N))$\,. In this case we can
use Kosaki's change of density argument (see \cite{Ko}) in order
to investigate the interpolation space $[I_{\infty,
p}^0({BMO}_c(\N), L_p^c{MO}(\N)]_\theta$\,, for $0< \theta< 1$\,.

\begin{prop}\label{rotz} Let $2\le p< q\le \infty$, $p\le v<\infty$ such that $\frac1p=\frac1q+\frac1v$ and $0<\theta<1$. Then for all $0\leq \eta< 1$\,, the
spaces $[\tilde{I}_{q,p}^{\eta}(L_q^cMO(\N)),
L_p^cMO(\N)]_{\theta}$ and $[\tilde{I}_{q,p}^1(L_q^cMO(\N)),
L_p^cMO(\N)]_{\theta}$ are isometrically isomorphic. The
isomorphism \[ T_{\theta}:[\tilde{I}_{q,p}^1(L_q^cMO(\N)),
L_p^cMO(\N))]_\theta\to [\tilde {I}_{q,p}^{\eta}(L_q^cMO(\N)),
L_p^cMO(\N)]_{\theta} \] is defined by
 \[ T_{\theta}((x_mD^{\frac1v})_m)\lel
  \left(D^{\frac{(1-\eta)(1-\theta)}{v}}x_mD^{\frac{1}{v}-\frac{(1-\eta)(1-\theta)}{v}}\right)_m\,,
  \]
for all adapted finite sequences $(x_m)_m\in L_q^cMO(\N)$\,.
\end{prop}

\begin{proof} We will first construct a contraction
 \[ T_{\theta}:[\tilde{I}_{q,p}^1(L_q^cMO(\N)), L_p^cMO(\N)]_{\theta}\to [\tilde{I}_{q, p}^{\eta}(L_q^cMO(\N)),L_p^cMO(\N)]_{\theta} \pl .\]
According to Lemma \ref{densarg}, it suffices to define
$T_{\theta}$ on the dense set $\bigcup\limits_{m \geq
0}\tilde{I}_{q,p}^1(\N_{m,a}D^{\frac1q})$ and show that the
restriction is contractive. Furthermore, by the density result
\cite[Lemma 4.3.2]{BL}, it is enough to consider elements of the
form
  \[ x \lel f(\theta)D^{\frac1p}\,, \]
where
 \[  f(z)\lel  \summ_{k=1}^l g_k(z) y_k \]
with $y_k\in \N_{m,a}$, for some $m$\,, and $g_k$ analytic
functions which vanish at infinity such that
 \[ \max\left\{\sup_{t\in \mathbb{R}}\noo f(it)D^{\frac1q}\rrm_{{L_q^cMO(\N)}},
  \sup_{t\in \mathbb{R}} \noo f(1+it)D^{\frac1p}\rrm_{L_p^cMO(\N)}\right\}\pl< \pl 1
 .\]
For all $z$ in the unit strip $S$ we define
 \begin{equation*}
 g(z)=
 D^{\frac{(1-\eta)(1-z)}{v}}f(z)D^{\frac1q+\frac{\eta}{v}+\frac{(1-\eta)z}{v}}
 =  \summ_{k=1}^l g_k(z)
 (D^{\frac{(1-\eta)(1-z)}{v}}y_kD^{\frac{(1-\eta)z}{v}}D^{
 \frac1q+\frac{\eta}{v}})
 \pl .
 \end{equation*}
According to Fact \ref{anfact}, this function is analytic, taking
values in $L_p(\N_m)\subset L_p^c(\N_m; l_\infty)$. Then, we
observe that $\si_t$ is an isometry on $L_p^cMO(\N)$ and therefore
 \begin{align*}
 \noo g(1+it)\rrm_{L_p^cMO(\N)}
 &= \noo
 D^{\frac{-it(1-\eta)}{v}}f(1+it)D^{\frac1q
 +\frac{\eta}{v}+\frac{1-\eta}{v}}
 D^{\frac{(1-\eta)it}{v}}\rrm_{L_p^cMO(\N)}
 \\
 &= \noo f(1+it)D^{\frac1p}\rrm_{L_p^cMO(\N)}< 1  \pl .
 \end{align*}
On the other hand, by (\ref{e.3.333}),
 \begin{align*}
 \noo D^{-\frac{1-\eta}{v}}g(it) D^{-\frac{\eta}{v}} \rrm_{{L_q^cMO(\N)}}
 &= \noo
 D^{\frac{-it(1-\eta)}{v}}
 f(it)D^{\frac{(1-\eta)it}{v}}
 D^{\frac1q}\rrm_{{L_q^cMO(\N)}}
 \\
 &= \noo
 D^{\frac{-it(1-\eta)}{v}}
 f(it)D^{\frac1q}D^{\frac{(1-\eta)it}{v}}
 \rrm_{{L_q^cMO(\N)}} \\
 &=
  \noo f(it)D^{\frac1q}\rrm_{{L_q^cMO(\N)}} < 1 \pl .
 \end{align*}
Therefore, for $x=f(\theta)D^{\frac1p}$ it follows that
 \begin{align}
  T_{\theta}(x)&=  g(\theta) \lel
  \summ_{k=1}^l g_k(\theta)
 D^{\frac{(1-\eta)(1-\theta)}{v}}y_kD^{\frac{(1-\eta)\theta}{v}}D^{\frac1q+\frac{\eta}{v}}
 \label{rist}\\
 &=
 D^{\frac{(1-\eta)(1-\theta)}{v}}f(\theta)D^{\frac1q}D^{\frac{(1-\eta)\theta}{v}+\frac{\eta}{v}}
  \lel
  D^{\frac{(1-\eta)(1-\theta)}{v}}f(\theta)D^{\frac1p-\frac{(1-\eta)(1-\theta)}{v}}
  \pl \nonumber
 \end{align}
is an element in the interpolation space
$[\tilde{I}_{q,p}^{\eta}(L_q^cMO(\N)), L_p^cMO(\N)]_{\theta}$\,.
Hence $T_{\theta}$ defined by (\ref{rist}) extends to a
contraction. Conversely, a similar argument shows that the map
$\tilde{T_\theta}$ defined by
\begin{align*}
  \tilde{T}_{\theta}(D^{\frac{1-\eta}{v}}f(\theta)D^{\frac1q}
 D^{\frac{\eta}{v}}) &=
  D^{\frac{(1-\eta)\theta}{v}}f(\theta)D^{\frac1q+\frac{\eta}{v}+\frac{(1-\eta)(1-\theta)}{v}
  }
 \end{align*}
is a contraction. Since $f$ is analytic, there exists $y\in
\N_{m,a}$ such that
 \[ D^{\frac{(1-\eta)\theta}{v}}f(\theta)\lel yD^{\frac{(1-\eta)\theta}{v}} \pl
 .\]
Furthermore, we deduce that
 \begin{align*}
 T_{\theta}\tilde{T}_{\theta}(D^{\frac{1-\eta}{v}}f(\theta)D^{\frac1q}
 D^{\frac{\eta}{v}})
 &= T_{\theta}(yD^{\frac1p}) \\
 &=
 D^{\frac{(1-\eta)(1-\theta)}{v}}yD^{\frac1p-\frac{(1-\eta)(1-\theta)}{v}}
 \\
 &=
 D^{\frac{(1-\eta)(1-\theta)}{v}}yD^{\frac{(1-\eta)\theta}{v}}D^{\frac1q+\frac{\eta}{v}}
 \\
 &=D^{\frac{1-\eta}{v}}f(\theta)D^{\frac1q+\frac{\eta}{v}} \pl.
 \end{align*}
Hence $T_{\theta}\tilde{T}_{\theta}=\text{id}$. Similarly, we may
show that $\tilde{T}_{\theta}T_{\theta}=\text{id}\,,$ and the
assertion follows.
\end{proof}

\begin{cor} \label{rotate} Let $2< p <s<q\le \8\,,
\,0<\theta < 1$ such that $ \frac1s\lel
\frac{1-\theta}{q}+\frac{\theta}{p}$, and $0\leq \eta\leq 1$\,.
Then, isometrically, \[ {[\tilde{I}_{q,p}^\eta(L_q^cMO(\N)),
L_p^cMO(\N)]_\theta} \lel
 {\tilde{I}_{s,p}^\eta}(L_s^cMO(\N))\,.\]
\end{cor}

\begin{proof} Let us note first that
$\frac1p-\frac1s=\frac{1-\theta}{v}$\,, where
$\frac1v={\frac1p}-{\frac1q}$\,. Then, we observe that
 \[ T_{\theta}(xD^{\frac1p-\frac1s})
 \lel D^{(1-\eta)\frac{1-\theta}{v}}x
 D^{\frac1p-\frac1s-(1-\eta)\frac{1-\theta}{v}}
 \lel D^{(1-\eta)(\frac{1}{p}-\frac{1}{s})}xD^{\eta(\frac{1}{p}-\frac{1}{s})}  \pl . \]
The assertion follows now from Proposition \ref{intlqcmo} and
Proposition \ref{rotz}.
\end{proof}

\noindent Recall now the row version of the martingale norm
(\ref{lpcmonorm})
 \[  \noo x\rrm_{L_p^rMO(\N)} \lel \sup_m \noo \sup_{n\le m}  \E_n\kla \summ_{k=n}^m d_k(x)d_k(x)^*\mer  \rrm_{\frac{p}{2}}^{\frac12}  \pl . \] \lz
Interchanging the roles of columns and rows we immediately get the
following result.\lz

\begin{prop} \label{rowv}
Let $2< p <t<q\le \8\,, 0< \theta <1$  such that
$\frac1s=\frac{1-\theta}{q}+{\frac{\theta}{p}}$, and $0\leq
\eta\leq 1$. Then, isometrically,
 \[
 {[\tilde{I}_{q,p}^\eta(L_q^rMO(\N)), L_p^rMO(\N)]_\theta}  =
 \tilde{I}_{s,p}^\eta(L_s^rMO(\N))\,. \]
\end{prop}

For $2< p\leq \infty$, the following martingale norms were defined
in \cite{JX}:
\[ \noo x\rrm_{L_pMO(\N)} \lel \max\{ \noo x \rrm_{L_p^cMO(\N)}, \noo
x\rrm_{L_p^rMO(\N)}\}\,, \] respectively,
\[ \|x\|_{H_p(\N)}=\max\{\|x\|_{H_p^c(\N)}\,,
\|x^*\|_{H_p^c(\N)}\}\,. \] We now state the first version of our
main result.

\begin{theorem}\label{int1}
Let $2< p <s<q\le \8\,, 0< \theta <1$ such that
$\frac1s=\frac{1-\theta}{q}+ {\frac{\theta}{p}}$, and $0\leq
\eta\leq 1$\,. Then, with equivalent norms,
 \[
 {[{I}_{q,p}^\eta(L_qMO(\N)), L_p(\N)]_\theta} =
 {I}_{s,p}^\eta(L_s(\N))\,. \]
\end{theorem}

\begin{proof} Recall that, by Corollary \ref{corapplic}, we have
\[ \|x\|_{L_p^c{MO}(\N)}\leq \lambda_p\|x\|_p\,, \]
where $\lambda_p=O(1)$ as $p\rightarrow \infty\,.$ Combining this
inequality with Corollary \ref{rotate}, Proposition \ref{rowv} and
Corollary 5.6 of \cite{Ra2}, it follows that
\begin{eqnarray*}
{[{I}_{q, p}^\eta(L_qMO(\N)), L_p(\N)]}_\theta &\subseteq
&{{[\tilde{I}_{q, p}^\eta(L_q^cMO), L_p^cMO]}_\theta}\cap
{{[\tilde{I}_{q,
p}^{\eta}(L_q^rMO), L_p^rMO]}_\theta}\\
&\subseteq &{\tilde{I}_{s, p}^\eta(L_s^cMO)}\cap {\tilde{I}_{s,
p}^\eta(L_s^rMO)}\\
&=& \tilde{I}_{s, p}^\eta({L_s^cMO}\cap {L_s^rMO})\\
&\subseteq &\tilde{I}_{s, p}^\eta (L_s(\N))\,,
\end{eqnarray*}
where the norm of the last inclusion above is
$\lambda^{\prime}_s\leq Cs$ as $s\rightarrow \infty\,;$ here $C$
is an absolute constant. The reverse inclusion follows by applying
Kosaki's interpolation result (\ref{Kosinterp})\,, together with
Corollary \ref{rotate} and Proposition \ref{rowv}, namely
\begin{eqnarray*}
I_{s, p}^\eta(L_s(\N))&\subseteq &{[I_{q, p}^\eta(L_q(\N)),
L_p(\N))]}_\theta\\&\subseteq &{[I_{q, p}^\eta(L_q^cMO\cap
{L_q^rMO}), L_p(\N)]}_\theta\\&=&{[I_{q, p}^\eta(L_qMO),
L_p(\N)]}_\theta\,.
\end{eqnarray*}\\[-1.4cm]
\end{proof}

As an application, we obtain  our main result of this section. The
following continuous inclusion is proved using the fact that we
have an interpolation scale of spaces.

\begin{theorem}\label{bmoinlp}  Let $1\leq p <\8$ and $0\leq \eta\leq 1$, then  the inclusion map
 \[ I_{\infty, p}^ {\eta}(BMO) \subset  L_p(\N) \]
is bounded, with norm $c(p)\leq cp$\,, where $c$ is an absolute
constant.
\end{theorem}\lz

\begin{proof} Assume that $p\geq 8$\,. Let $4< s< {p/2} <\8$ and
define $0<\theta<1$ such that
$\frac1p=\frac{1-\theta}{\8}+\frac{\theta}{s}$. Since the complex
interpolation method is an exact interpolation functor of exponent
$\theta$, we deduce from (the proof of) Theorem \ref{int1}\,, that
for a finite martingale difference sequence $x=(d_k)_{k=0}^n \in
BMO$ we have
\begin{eqnarray*}
\|D^{\frac{1-\eta}{p}}xD^{\frac{\eta}{p}}\|_p&\leq
&c(p)\|D^{\frac{1-\eta}{s}}xD^{\frac{\eta}{s}}\|_{[I_{\infty,
s}^{\eta}(BMO),
 L_{s}(\N)]_\theta}\\&\leq& c(p)\|x\|_{BMO}^{1-\theta}\|D^{\frac{1-\eta}{s}}xD^{\frac{\eta}{s}}\|_{s}^{\theta}\leq
 c(p)\|x\|_{BMO}^{1-\theta}\|D^{\frac{1-\eta}{p}}xD^{\frac{\eta}{p}}\|_p^{\theta}\,,
 \end{eqnarray*}
where $c(p)\leq Cp$, with $C$ an absolute constant. Thus, we get
\begin{equation}\label{eq:76756}
\|D^{\frac{1-\eta}{p}}xD^{\frac{\eta}{p}}\|_p\leq
(c(p))^{\frac1{1-\theta}}\|x\|_{BMO}\,.
\end{equation}
Under the hypotheses on $p$ and $s$\,, a simple computation shows
that
\[ (c(p))^{\frac{1}{1-\theta}}\leq cp\,, \] where $c$ is an absolute constant. Therefore (\ref{eq:76756})
implies that
\[ \|D^{\frac{1-\eta}{p}}xD^{\frac{\eta}{p}}\|_p\leq cp\ \|x\|_{BMO}\,. \]
For an arbitrary element $x=(d_k)_{k\geq 0}$ in $BMO$, we consider
the sequence of projections $P_n(x)\,,$ where
\[ P_n(x)=(d_k)_{k=0}^n\,, \]
whose images $I_{\infty, p}^{\eta}(P_n(x))$ are uniformly bounded
and thus converge in $L_p(\N)$ to $I_{\infty, p}^{\eta}(x)$. This
yields the assertion for $p\geq 8$\,. Since $L_8(\N)$ embeds into
$L_p(\N)$ for all $1\leq p< 8$\,, the proof is complete.
\end{proof}

Note that the symmetric embedding
\[ \Psi_{p}(x) \lel D^{\frac{1}{2p}}xD^{\frac{1}{2p}} \]
is positivity preserving and thus, it can be considered as the
natural embedding of $\N$ into $L_p(\N)$. As a special case of
Theorem \ref{bmoinlp} we obtain the following

\begin{cor}\label{symbmoinlp} If $1\leq p< \infty$, then the
inclusion map
 \[ \Psi_{p}(BMO) \subset  L_p(\N) \]
is bounded, with norm $c(p)\leq cp$, where $c$ is an absolute
constant.
\end{cor}


\begin{rem}\label{rem95}
{\rm The following example which answers a question raised by Tao
Mei, shows that the above inclusion does not hold if we only
restrict to column, respectively row versions of the ${BMO}$ norm.
 Let $n$ be
a positive integer and consider the von Neumann algebra
\[ \N={L_\infty([0,1])}\bar{\otimes}M_n\,, \]
where $M_n$ is the algebra of $n\times n$ complex matrices. For
$k\geq 1$ let $\Sigma_k$ be the $\sigma$-algebra generated by
dyadic intervals in $[0, 1]$ of length $2^{-k}$\,. Denote by
$\N_k$ the subalgebra ${L_\infty([0, 1],
\Sigma_k)}\bar{\otimes}M_n$ of $\N$ and let
$\E_k={\mathbb{E}_k}\otimes {\text{Id}_{M_n}}$ be the conditional
expectation onto $\N_k$. Let $\varepsilon_k\in\{-1, 1\}$ and
define
\[ x=\sum\limits_{k=1}^n \varepsilon_k \otimes e_{1k}\,. \]
Then $x$ is a martingale relative to the filtration $(\N_k)$\,.
The martingale differences are given by $d_k(x)=
\varepsilon_k\otimes e_{1k}$\,. A simple computation shows that
$\|x\|_{L_p(\N)}=n^{\frac12-\frac1p}$\,, while
\[ \|x\|_{{BMO}_c}=\sup\limits_m \left\|\sum\limits_{k=m}^n
\E_m(d_k^*(x)d_k(x))\right\|_\infty^\frac12=\sup\limits_k\|\varepsilon_k^2\|_\infty^\frac12=1\,.
\] Assume that $p>
2$\,. Then for any $C> 0$\,, there exists $n\geq 1$ such that
$n^{\frac12-\frac1p}> C$\,. This implies that ${BMO}_c(\N)$ is not
contained in $L_p(\N)$\,.
 } \end{rem}

\section{Main results}

In this section, we will show how to derive the John-Nirenberg
type result from Theorem \ref{bmoinlp}. Let us start with an
immediate application.

\begin{cor}\label{reform}  Let $1\leq p< \infty$\,, a positive integer $n$ and an element $a\in L_p(\N_n)$\,.
Then, there exists an absolute constant $c>0$ such that for all
$x\in \N$\,,
\begin{equation}\label{inclusbmolp}
\noo (x-x_{n-1})a\rrm_p \kl c p \noo x\rrm_{BMO} \noo a\rrm_p\,.
\end{equation}
\end{cor}

\begin{proof} Assume that $p\geq 2$\,. It
suffices to show (\ref{inclusbmolp}) for positive elements $a\in
L_p(\N_n)$ with $\|a\|_p=1$\,. Furthermore, by approximation with
elements of the form $(a^p+ \eps D)^{\frac1p}$, we may assume that
$a$ has full support. Define a new state $\phi_a$ on $\N$ by
\[ \phi_a(x)\lel \text{tr}(a^px)\,, \quad x\in \N\,. \]
Denote by $L_p(\N, \phi_a)$ the noncommutative $L_p$ space
associated to the state $\phi_a\,.$ Note that
\begin{equation}\label{lpfia}
\|x-x_{n-1}\|_{L_p(\N, \phi_a)}=\|(x-x_{n-1})a\|_p\,.
\end{equation}
For all non-negative integers $k$ define
$\widetilde{\N}_k=\N_{n+k}$\,. Then $\widetilde{\N}_0=\N_n$ and
$(\widetilde{\N}_k)_{k\ge 0}$ is a filtration of $\N$\,. Moreover,
if $m\ge n$ and $x\in \N$\,, then $xa^p\in L_1(\N_m)$\,. Hence,
\begin{equation}\label{newstate}
 \phi_a(yx)\lel \text{tr}(a^pyx) \lel \text{tr}(yxa^p) \lel \text{tr}(\E_m(y)xa^p)
 \lel \phi_a(\E_m(y)x) \pl ,
 \end{equation}
 for all $y\in \N_m$\,. Denote by
$BMO((\widetilde{\N}_k)_{k\geq n}\,, \phi_a)$ the $BMO$ space
associated to the filtration $(\widetilde{\N}_k)_{k\geq 0}$ of
$\N$ and the state $\phi_a$\,. From Theorem \ref{bmoinlp} it
follows that
\begin{equation}\label{eq89889}
\noo x-x_{n-1}\rrm_{L_p(\N, \phi_a)} \kl c p \noo
 x-x_{n-1}\rrm_{BMO((\widetilde{\N}_k)_{k\geq n}\,, \phi_a)}\,,
 \end{equation}
The modular automorphism group of $\phi_a$ is given by
\[ \si_t^{\phi_a}(x)\lel a^{ipt}xa^{-ipt}\,, \quad t\in \mathbb{R}\,. \]
From (\ref{newstate}) we deduce that for all $k\geq 0$\,, the
conditional expectations $\tilde{\E_k}=\E (\cdot|
\widetilde{\N}_k)$ associated with the state $\phi_a$ are given by
the original conditional expectations, namely, $\tilde{\E}_k\lel
\E_{k+n}$\,. Therefore,
 \[ \tilde{\E}_k(x-x_{n-1})\lel  \E_{k+n}(x-x_{n-1}) \lel x_{k+n}-x_{n-1}\,. \]
Hence, denoting $x-x_{n-1}$ by $y$\,, we obtain
 \begin{align*}
\tilde{\E}_k((y-\tilde{\E}_{k-1}(y))^*(y-\tilde{\E}_{k-1}(y)))=\E_{n+k}((x-x_{n+k-1})^*(x-x_{n+k-1}))\,.
 \end{align*}
By the definition of the $BMO$ norm, it follows that
 \[ \noo x-x_{n-1}\rrm_{BMO((\widetilde{\N}_k)_{k\geq n}\,, \phi_a)}\kl \noo x\rrm_{BMO} \pl. \]
Together with (\ref{eq89889}) and (\ref{lpfia}), this yields the
assertion in the case $p\geq 2$\,. Furthermore, if $1\leq p< 2$
and $0\leq a\in L_p(\N)$ with $\|a\|_p=1$\,, consider
$a_1=a^{\frac12}$\,. Then $0\leq a_1\in L_{2p}(\N)$ and
$\|a_1\|_{2p}=1$\,. By H\"{o}lder's inequality it follows that
\begin{equation}\label{hol}
\|(x-x_{n-1})a\|_p\leq \|(x-x_{n-1})a_1\|_{2p}\,.
\end{equation}
Note that $2p\geq 2$\,. Applying (\ref{inclusbmolp}) and
(\ref{hol}) to $a_1$, we obtain the conclusion.
 \end{proof}

\begin{rem}\label{intconditionlp}
{\rm Note that for any integer $n\geq 0$\,, the inclusion
$\N\subset L_\infty^c(\N,\E_n)$ is injective. Thus, given $2\leq
p< \infty$\,, we can consider the space
$[\N,L_\infty^c(\N,\E_n)]^{\frac{2}{p}}$ obtained by the upper
method of complex interpolation of exponent $\frac{2}{p}$\,. These
are particular examples of conditional $L_p$-spaces considered in
\cite{JPa}\,. Following ideas of Pisier \cite{Pi6}, it is proved
in \cite{JPa} that
\begin{equation}\label{interpnormjpa}
\|x\|_{[\N,L_\infty^c(\N,\E_n)]^{\frac{2}{p}}}=\sup_{\noo
a\rrm_{L_p(\N_n)}\le 1} \noo xa\rrm_{L_p(\N)}\,. \end{equation}
Therefore, as already mentioned in the introduction, the
${{BMO}_p^c}$ norms defined by (\ref{eq4534}) are in fact
interpolation norms, namely
\begin{equation}\label{eqnormbmopc}
\|x\|_{{BMO}_p^c}=\|x-x_{n-1}\|_{[\N,L_\infty^c(\N,\E_n)]^{\frac{2}{p}}}\,.
\end{equation}
 }
\end{rem}

\begin{proof}[{\it Proof of Theorem \ref{martingal}}] Let $2< p< \infty$\,. Trivially, we have
 \[ {[\N,L_\infty^c(\N, \E_n)]}^{\frac{2}{p}}\subset
 L_\infty^c(\N,\E_n) \pl. \]
Therefore, we deduce that
 \[ \noo \E_n((x-x_{n-1})^*(x-x_{n-1}))\rrm_\infty^{\frac12}\kl \noo
 x-x_{n-1}\rrm_{[\N, L_\infty^c(\N,\E_n)]^{\frac{2}{p}}} \pl .\]
Combining this with (\ref{eqnormbmopc})\,, it follows immediately
that
 \[ \noo x\rrm_{BMO}\kl \noo x\rrm_{BMO_p} \pl .\]
For the converse, we deduce from (\ref{eq4534}) and Corollary
\ref{reform} that
 \[ \|x\|_{{BMO}_p^c}
 \lel \sup_{\noo a\rrm_{L_p(\N_n)}\le 1} \noo (x-x_{n-1})a\rrm_p
 \kl c p \noo x\rrm_{BMO} \pl .\]
Applying  the same argument to $x^*$, we deduce that
 \[ \noo x\rrm_{BMO_p}\kl c p \noo x\rrm_{BMO} \pl ,\]
which concludes the proof.
\end{proof}

\begin{rem}\label{lexp}
{\rm As an application of Theorem \ref{martingal} we now prove
that the noncommutative analogue of the classical result
\begin{equation}\label{bmointolexp}
BMO\subset L_\text{exp}
\end{equation}
holds in the setting of a semifinite von Neumann algebra $\N$
equipped with a normal, faithful trace $\tau$\,. In this context
the space $L_\text{exp}(\N)$ can be defined following the general
scheme of {\em symmetric spaces} associated to $(\N, \tau)$ and a
rearrangement invariant  Banach function space developed in
\cite{DDP1} and \cite {DDP2}. Therefore, motivated by the
classical definition of the Zygmund space $L_{\text{exp}}$\,, we
define
\begin{equation}\label{normlexp}
\|x\|_{L_{\text{exp}}(\N)}=\inf\{\lambda> 0:
\tau(e^{{\frac{|x|}{\lambda}}-1})\leq 1\}\,.
\end{equation}
Suppose that $\|x\|_{BMO}\leq 1$ and let $\lambda> 0$\,. Then,
using the power series expansion for the exponential, it follows
from Theorem \ref{martingal} that
\begin{eqnarray*}
\tau(e^{{\frac{|x|}{\lambda}}-1}) =\frac1e\sum\limits_{k=0}^\infty
\frac{\tau(|x|^k)}{{\lambda^k}{k!}}&=&\frac1e\left(1+\sum\limits_{k=1}^\infty
\frac{\|x\|_k^k}{{\lambda^k}{k!}}\right)\\&\leq&
\frac1e\left(1+\sum\limits_{k=1}^\infty
\frac{{c^k}{k^k}}{{\lambda^k}{k!}}\right)\\&\leq&
\frac1e\left(1+\sum\limits_{k=1}^\infty
\frac{{c^k}{e^k}}{\lambda^k}\right)=\frac1e
\sum\limits_{k=0}^\infty \left(\frac{ce}{\lambda}\right)^k\,.
\end{eqnarray*}
For the last inequality we have used the fact, which can be
verified by induction, that $k^k\leq {{k!}{e^k}}$\,, for all
integers $k\geq 1$\,. A simple computation shows that if
$\frac1\lambda\leq \frac1{ce}\left(1-\frac1e\right)$\,, then
$\frac1e \sum\limits_{k=0}^\infty
\left(\frac{ce}{\lambda}\right)^k\leq 1$\,. By (\ref{normlexp})\,,
it follows that
\[ \|x\|_{L_{\text{exp}}(\N, \tau)}\leq
\frac1{ce}\left(1-\frac1e\right)\,. \] Hence $x\in
{L_{\text{exp}}(\N)}$ and (\ref{bmointolexp}) is proved. }
\end{rem}

\noindent {\bf Open problem}: Is there a universal constant $C$
such that if $x\in {BMO}$\,, then
\[ \sup\limits_n \|\E_n(|x-x_{n-1}|^{p})\|^{\frac1p}_\infty\leq Cp\|x\|_{BMO}\,, \]
for all $1\leq p< \infty$\,?

The noncommutative versions of function space $BMO$ have been
developed by Mei
 \cite{Mei}. Namely, for a function $x:\rz \to \N$\,, define the following norm
 \[ \noo x\rrm_{{BMO}_c} \lel \sup_{\emptyset \neq I} \noo \intt_{I}
 (x(t)-x_I)^*(x(t)-x_I) \frac{dt}{|I|} \rrm_{\N} \pl ,\]
 where $x_I=\intt_{I} x(t)\frac{dt}{|I|}$\,, for every nonempty
 interval $I$\,.
Furthermore, define
 \[ \noo x\rrm_{{BMO}_r} \lel \noo x^*\rrm_{{BMO}_c} \]
and
 \[ \noo x\rrm_{BMO} \lel \max\{\noo x\rrm_{{BMO}_c},\noo
 x\rrm_{{BMO}_r}\} \pl .\]

As in the martingale setting, we use interpolation to define the
$BMO_p$-norms in this context. Namely, for $2< p< \infty$ and a
fixed nonempty interval $I$\,, let $\mu_I=\frac{dt}{|I|}$ and
define the space
 \[ \N(L_p^c(\mu_I))\lel
 [\N\bar{\ten}L_\infty(\mu_I), \N\bar{\ten}L_2^c(\mu_I)]^{\frac{2}{p}}
 \pl \]
by (the upper method of) complex interpolation. Recall that
$\N\bar{\ten}L_\infty(\mu_I)=L_{\infty}(I,\mu_I; \N)$ is the usual
space of measurable, essentially bounded $\N$-valued functions,
while the space $\N\bar{\ten}L_2^c(\mu_I)$ is defined following
the operator space tradition of viewing $L_2(\mu_I)$ implemented
as a column in $B(L_2(\mu_I))$. Equivalently,
 \[ \noo x\rrm_{\N\bar{\ten}L_2^c(\mu_I)} \lel \noo \intt_I
 (x(t))^*x(t){d\mu_I(t)}\rrm_{\N}^{\frac12} \pl \]
for all $\N$-valued functions which are measurable with respect to
the strong operator topology. Note that $\E_I(x)\lel \int_I x(t)
d\mu_I(t)$
is a conditional expectation onto
$\N_I=\N\bar{\ten}L_\infty(\mu_I)$\,. Therefore, the spaces
$\N(L_p^c(\mu_I))$ are special cases of conditional $L_p$-spaces
considered in \cite{JPa}. In particular the following norm
estimate holds
\begin{equation} \label{form}
  \noo 1_Ix\rrm_{\N(L_p^c(\mu_I))} \lel \sup_{\noo a\rrm_{L_p(\N_I)}\le
 1} \kla \intt_{I} \noo x(t)a\rrm_{L_p(\N)}^p {d\mu_I(t)}
 \mer^{\frac1p} \pl .\end{equation}

As in \cite{Mei}, we will use the following result to transfer
martingale results to intervals.

\begin{lemma}[T. Mei]\label{mei} Let $I\subset \rz$ be an interval. Then there exists
a dyadic interval $J=[\frac{k-1}{2^n},\frac{k}{2^n}]$ or
$J'=[\frac{1}{3\cdot 2^{n'}}+\frac{k'-1}{2^{n'}},\frac{1}{3\cdot
2^{n'}}+\frac{k'}{2^{n'}}]$ such that
 \[ I\subset J \quad \mbox{and} \quad |J|\kl 6|I| \]
or
 \[ I\subset J' \quad \mbox{and}\quad  |J'|\kl 6|I| \pl .\]
\end{lemma}

\begin{lemma}\label{elem1} For a nonempty interval $I\subset \rz$ and an element  $x\in \N\bar{\ten}L_\infty(\rz)$ let
 \[ \noo x\rrm_{p,I}\lel \sup_{\noo a\rrm_{L_p(\N)}\le 1} \noo
 (1_Ix-x_I)a\rrm_{L_p(I,\mu_I;L_p(\N))} \pl .\]
Then, for all intervals $I\subset J \subset \rz$\,,
 \[ \noo x\rrm_{p,I}\kl 2 \kla \frac{|J|}{|I|}\mer^{\frac1p} \noo
 x\rrm_{p,J} \pl .\]
\end{lemma}

\begin{proof} This is of course a standard argument. Let $a\in
L_p(\N)$. Then we have
 \begin{align*}
  \noo (x_I-x_J)a\rrm_p &= \noo \intt_I [(x(s)-x_J)a] \,{d\mu_I(s)}
  \rrm_p \kl
   \kla \intt_I \noo (x(t)-x_J)a\rrm_p^p {d\mu_I(t)}
  \mer^{\frac1p} \pl .
  \end{align*}
Therefore, we obtain
 \begin{align*}
 \noo (1_Ix-x_I)a\rrm_{L_p(I,\mu_I;L_p(\N))} &\le \noo (1_Ix-x_J)a\rrm_{L_p(I,\mu_I;L_p(\N))}+
  \noo (x_I-x_J)a\rrm_p \\
  &\le 2 \kla \intt_{I}\noo (x(t)-x_J)a\rrm_p^p {d\mu_I(t)}
  \mer^{\frac1p}\\
  &\le 2 \kla \frac{|J|}{|I|}\mer^{\frac1p}
  \kla \intt_{J}\noo (x(t)-x_J)a\rrm_p^p {d\mu_J(t)}
  \mer^{\frac1p} \pl .
  \end{align*}
Taking supremum over all elements $a$ in the unit ball of
$L_p(\N)$ implies the assertion.
\end{proof}

We are now ready to formulate a noncommutative interval version of
the John-Nirenberg theorem.

\begin{theorem}\label{interval} For all $2<p<\infty$, there exists an absolute constant
$c>0$ such that
 \[ \noo x\rrm_{BMO}\kl \sup_{\emptyset \neq I\ \text{interval}} \ \max\{\noo
 x-x_I\rrm_{\N(L_p^c(\mu_I))},
 \noo
 x^*-x^*_I\rrm_{\N(L_p^c(\mu_I))}\} \kl cp \noo x\rrm_{BMO} \pl
 .\]
\end{theorem}

\begin{proof}
Denote by $\widetilde{\N}$ the von Neumann algebra tensor product
$\N\bar{\ten} L_{\infty}[0,1]$. For each integer $n\geq 0$, let
$\Si_n$ be the $\si$-algebra generated by dyadic intervals in $[0,
1]$ of length $2^{-n}$\,. Denote by $\widetilde{\N}_n$ the
subalgebra $\N\bar{\ten} L_{\infty}([0,1],\Si_n)$ of
$\widetilde{\N}$ and let $\tilde{\E}_n$ be the conditional
expectation onto $\widetilde{\N}_n$\,. Consider the $BMO$ space
associated to the filtration $(\widetilde{\N}_n)_{n\geq 0}$ of
$\widetilde{\N}\,.$ For $x\cong (x_n)_{n\geq 0}\,,$ where
$x_n=\tilde{\E}_n(x)$\,, we deduce that
 \[  \sup_n \ \max\{ \noo (x-x_n)\rrm_{L_p^c(\tilde{\N}, \tilde{\E}_n)},
 \noo (x-x_n)^*\rrm_{L_p^c(\tilde{\N},\tilde{\E}_n)}\} \kl \noo x\rrm_{BMO}
 \pl . \]
An application of (\ref{form}) shows that
 \begin{align*}
  \noo x-x_n\rrm_{L_p^c(\tilde{\N}, \tilde{\E}_n)}^p &=
  \sup_{\summ_k \noo a_k\rrm_p^p |I_k|\le 1}
  \summ_{k} |I_k| \intt_{I_k}\noo (x(t)-x_{I_k})a_k\rrm_p^p {d\mu_{I_k}(t)}\\
  &= \sup_k \noo 1_{I_k}x-x_{I_k}\rrm_{\N(L_p^c(\mu_{I_k}))}^p \pl .
  \end{align*}
Here the supremum is taken over all intervals $I_k\subset [0,1]$
of length $2^{-n}$. Therefore, we may reformulate the martingale
result in Theorem \ref{martingal} as
 \begin{align*}
 &  \sup_{I \mbox{ \scriptsize dyadic }\subset [0,1] } \max\{\noo 1_Ix-x_I\rrm_{\N(L_p^c(\mu_I))},
 \noo (1_Ix-x_I)^*\rrm_{\N(L_p^c(\mu_I))}\}\\
 &\pll \le
  cp \sup_{I \mbox{ \scriptsize dyadic }\subset [0,1] } \max\{\noo 1_Ix-x_I\rrm_{\N(L_2^c(\mu_I))},
 \noo (1_Ix-x_I)^*\rrm_{\N(L_2^c(\mu_I))}\} \pl .
 \end{align*}
Using the dilation and translation invariance properties of the
Lebesgue measure, it follows for all integers $k\geq 0$ that
 \begin{align*}
 &  \sup_{I \mbox{ \scriptsize dyadic }\subset [2^{-k},2^k] } \max\{\noo 1_Ix-x_I\rrm_{\N(L_p^c(\mu_I))},
 \noo (1_Ix-x_I)^*\rrm_{\N(L_p^c(\mu_I))}\}\\
 &\pll \le
  cp \sup_{I \mbox{ \scriptsize dyadic }\subset [2^{-k},2^k] } \max\{\noo 1_Ix-x_I\rrm_{\N(L_2^c(\mu_I))},
 \noo (1_Ix-x_I)^*\rrm_{\N(L_2^c(\mu_I))}\} \pl .
 \end{align*}
Taking the supremum over $k$, we obtain the desired inequality for
arbitrary dyadic martingales. Applying Mei's result \ref{mei}\,,
together with Lemma \ref{elem1} the assertion follows.
\end{proof}

We now discuss an analogue of the classical large deviation
inequality (\ref{eqlargedi}) in this setting. Our main tool is the
following corollary to a noncommutative version of Chebychev
inequality proved by Defant and Junge (see \cite{DJ}).

\begin{cor}[Defant-Junge]\label{andreasma}
Let $2\leq p< \infty$ and $x\in L_p(\N)$\,. Then, for every
$\varepsilon> 0$, there exists a projection $f\in \N$ with
$\phi(1-f)\leq \varepsilon$, such that whenever $yD^{\frac1p}=wx$,
where $w, y\in \N$\,, then
\begin{equation}\label{chebine}
\sup\limits_n \|f\E_n(y^*y)f\|^{\frac12}\leq (1+\sqrt{10})
\varepsilon^{-\frac1p}\|w\|_\infty\|x\|_p\,.
\end{equation}
Moreover,
\begin{equation}\label{chebineq}
\sup\limits_n\|\E_n(y)f\|\leq (1+\sqrt{10})
\varepsilon^{-\frac1p}\|w\|_\infty\|x\|_p\,.
\end{equation}
\end{cor}

\begin{proof}[{\it Proof of Theorem \ref{noncomldiv}}]
Recall that the state $\phi$ satisfies condition
(\ref{e.1.45679})\,. Let $c_2=e$\,. Furthermore, consider $t\geq
\frac1{c_1}$ and let $\varepsilon:=e^{-tc_1}$\,, where $c_1$ is a
constant to be made precise later. Define $p=4\ln
(\varepsilon^{-1})=4tc_2$\,. By construction it follows that
$4\leq p< \infty$\,. Now let $x\in {BMO}$ with $\|x\|_{BMO}\leq
1\,.$ Denote by $y$ the element $(x-x_0)D^{\frac1p}$\,. Note that
$y\in L_p(\N)$\,. By the noncommutative Doob's inequality (see
\cite{Ju}) and Theorem \ref{martingal}, we get
\begin{equation}\label{eq:something}
\|\sup\limits_n \E_n(y^*y)\|_{\frac{p}{2}}\leq
d_{\frac{p}{2}}\|y\|^2_p \leq d_{\frac{p}{2}} \,c^2 p^2\
\|x-x_0\|_{BMO}^2\,,
\end{equation}
where $c$ is the universal constant in Theorem \ref{martingal} and
$d_{\frac{p}{2}}$ is the constant in the noncommutative Doob's
inequality. Since $\frac{p}{2}\geq 2$\,, it follows that
$d_{\frac{p}{2}}=2$ (see \cite{Ju})\,. Furthermore, note that
\begin{equation}\label{bmoxnorm}
\|x-x_0\|_{BMO}\leq \|x\|_{BMO}\,.
\end{equation}
Indeed, we have the estimates
\[
\E_0((x-x_0)^*(x-x_0))=\E_0(x^*x-x^*x_0-\E_0(x^*)x+\E_0(x^*)\E_0(x))=\E_0(x^*x)-\E_0(x^*)\E_0(x)\,.
\]
By Kadison's inequality it follows that
\begin{equation}\label{eq:7828}
\|\E_0((x-x_0)^*(x-x_0))\|\leq \|\E_0(x^*x)\|\,.
\end{equation}
Moreover, for $n\geq 1$ we have
\[ (x-\E_0(x))-\E_{n-1}(x-\E_0(x))=(x-\E_0(x))-(\E_{n-1}(x)-\E_0(x))=x-\E_{n-1}(x)\,.
\]
Combining this with (\ref{eq:7828}) we deduce that
\[ \|x-x_0\|_{{BMO}_c}\leq \|x\|_{{BMO}_c}\,, \]
which yields the inequality (\ref{bmoxnorm}).\\
A further application of Kadison's inequality and
(\ref{eq:something}), together with (\ref{bmoxnorm}) implies that
\[ \|\sup\limits_m \E_m(y^*)\E_m(y)\|_{\frac{p}{2}}\leq 2c^2p^2\,. \]
Then by \cite{Ju}\,, Remark 3.7, there exists $0\leq a\in L_p(\N)$
and a sequence of positive contractions $w_m\in \N$ such that
\[ \E_m(y^*)\E_m(y)=aw_ma\,, \qquad
\|a\|_p^2\sup\limits_m\|w_m\|_\infty=\|\sup\limits_m
\E_m(y^*)\E_m(y)\|_{\frac{p}{2}}\leq \sqrt{2}cp\,. \] Using the
polar decomposition of $\E_m(y)$\,, we find partial isometries
$u_m$ such that, denoting $u_mw_m^{\frac12}$ by $z_m$\,, we obtain
\begin{equation}\label{estim}
\E_m(y)=z_m a\,, \qquad \|a\|_p\sup\limits_m\|z_m\|_\infty\leq
\sqrt{K}p\,.
\end{equation}
Note that  $\E_m(x)-x_0\in \N\,,$ for all non-negative integers
$m$\,. Indeed, for all $n\geq 0$\,,
\[ 1\geq \|x\|_{{BMO}_c}^2=\sup\limits_n \left\|\sum\limits_{k\geq n}
\E_n(d_k^* d_k)\right\|_\infty \geq
\|\E_n(d_n^*d_n)\|_\infty=\|d_n^*d_n\|_\infty\,.
\]
This implies that $d_n^*d_n\in \N$ and hence $d_n\in \N$\,.
Therefore $\E_m(x)-x_0=\sum\limits_{k=1}^m d_k\in \N\,.$ Since
\[ \E_m(y)=(\E_m(x)-x_0)D^{\frac1p}\,, \]
we are now in the condition of applying Corollary \ref{andreasma}
and deduce the existence of a projection $f\in \N$ with
$\phi(1-f)\leq \varepsilon$ such that for all $m\geq 0$\,,
\[ \|(\E_m(x)-x_0)f\|\leq
(1+\sqrt{10})\varepsilon^{-\frac1p}\|z_m\|_\infty\|a\|_p\,. \]
Using the estimates (\ref{estim}) we deduce that
\[ \|(\E_m(x)-x_0)f\|\leq Kp\ \varepsilon^{-\frac1p}\,, \]
where $K=(1+\sqrt{10})\sqrt{2}c\,.$ Therefore
$\|(\E_m(x)-x_0)h\|\leq Kp\ \varepsilon^{-\frac1p}\,,$ for all
$h\in fH$\,. Since the sequence $(\E_m)_{m\geq 0}$ converges to
$\text{Id}_\N$ in the strong operator topology, it follows that
\begin{equation}\label{estim2}
\|(x-x_0)f\|\leq Kp\ \varepsilon^{-\frac1p}\,. \end{equation}
 We now make precise the constant $c_1$\,. Namely, let
\[
c_1:=\frac1{4Ke^{1/4}}=\frac1{4\sqrt{2}(1+\sqrt{10})ce^{1/4}}\,.\]
 A simple computation yields the
 equality $K p\varepsilon^{-\frac1p}=t$\,. By (\ref{estim2}) we
deduce the inequality $\|(x-x_0)f\|\leq t\,.$ Moreover,
$\phi(1-f)\leq
\varepsilon< c_2e^{-tc_1}$\,, which yields the assertion.\\
For $0< t< \frac1{c_1}$\,, let $f=0$\,. The condition
$\|(x-x_0)f\|\leq t$ is then automatically satisfied. Furthermore,
note that
\[ \phi(1-f)=\phi(1)=1< \frac{e}{e^{tc_1}}=c_2e^{-tc_1}\,. \]
This completes the proof.
\end{proof}

\begin{rem}
{\rm After completing this paper, we have been informed that Tao
Mei has obtained a simple proof of the inequalities (\ref{JN9}),
but only in the tracial setting, based on the interpolation result
proved in \cite{Mu}\,. }
\end{rem}

\vspace*{0.4cm}


\begin{thebibliography}{456}

\bibitem{Ba}{\sc R.F. Bass}: {\em Probabilistic tecniques in
analysis}, Springer-Verlag, New York 1995.
\bibitem{Bu}{\sc D. Burkholder}: {\em Martingales and singular integrals in Banach spaces}, Handbook on the Geometry of Banach Spaces, Volume 1, W.B. Johnson and J. Lindenstrauss editors, North Holland (2001), 233-269.
\bibitem{BL}{\sc J. Bergh, J. L\"{o}fstr\"{o}m}: {\em Interpolation Spaces. An Introduction}, Springer-Verlag, New-York 1976.
\bibitem{CK}{\sc E. A. Carlen and P. Kr\'{e}e}, {\em On martingale inequalities in noncommutative stochastic analysis}, J. Funct. Analysis {\bf 158} (1998), 475-508.
\bibitem{BeSh}{\sc C. Bennett, R. Sharpley}: {\em Interpolation of operators}, Academic Press 1988.
\bibitem{DDP1}{\sc P.G. Dodds, T.K. Dodds, B. de Pagter}: {\em
Non-commutative Banach function spaces}, Math. Z. {\bf 201}
(1989), 583-597.
\bibitem{DDP2}{\sc P.G. Dodds, T.K. Dodds, B. de Pagter}: {\em Fully symmetric operator spaces}, Integr. Eq. Op. Theory Vol {\bf 15} (1992), 942-972.
\bibitem{DJ}{\sc A. Defant, M. Junge}: {\em Riesz summation in noncommutative $L_p$-spaces}, Preprint.
\bibitem{GR}{\sc J. Garcia-Cuerva, J. Rubio de Francia}: {\em
Weighted norm inequalities and related topics}, North Holland
Mathhematics Studies {\bf 116}, Notas de Matematica 104, Elsevier
Sci. Publishers B.V., 1985.
\bibitem{Gar}{\sc J.B. Garnett}: {\em Bounded analytic functions},
Academic Press, Inc. (London) Ltd., 1981.
\bibitem{Ga}{\sc A. Garsia}: {\em Martingale Inequalities: Seminar Notes on Recent Progress}, Math. Lecture Notes Series, 1973.
\bibitem{Gous}{\sc N. Ghoussoub, G. Godefroy, B. Maurey, W.
Schachermayer}: {\em Some topological and geometric structures in
Banach spaces}, Memoirs Amer. Math. Soc. {\bf 379} (1987).
\bibitem{Jo}{\sc F. John, L. Nirenberg}: {\em On functions of
bounded mean oscillation}, Comm. Pure Appl. Math. {\bf 167}
(1961), 416-426.
\bibitem{Ju}{\sc M. Junge}, {\em Doob's inequality for non-commutative martingales},
J. reine angew. Math. {\bf 549} (2002), 149-190.
\bibitem{Ju1}{\sc M. Junge},  {\em Fubini's theorem for ultraproducts of noncommutative $L_p$ spaces}; to appear in Canadian J. Math.
\bibitem{JPa}{\sc M. Junge, J. Parcet}, {\em Interpolation of amalgamated non-commutative $L_p$
spaces}, Preprint.
\bibitem{JS}{\sc M. Junge, D. Sherman}, {\em Noncommutative $L_p$
modules}, To appear in J. Operator Theory.
\bibitem{JX}{\sc M. Junge and Q. Xu}, {\em Non-commutative Burkholder/Rosenthal inequalities}, Ann.
Probab. {\bf 31} (2003) no. 2, 948-995.
\bibitem{JX2}{\sc M. Junge and Q. Xu}, {\em Non-commutative Burkholder/Rosenthal inequalities: Applications}, In preparation.
\bibitem{JX1}{\sc M. Junge, Q. Xu}: {\em The optimal orders of growth of the best constants in some non-commutative martingale inequalities}, Preprint, 2001.
\bibitem{KR}{\sc R.V. Kadison, J.R. Ringrose}: {\em Fundamentals of the Theory of Operator Algebras I, II}, Academic Press, 1986.
\bibitem{Koo}{\sc P. Koosis}: {\sc Introduction to $H_p$-spaces},
Cambridge tracts in mathematics {\bf 115}, Cambridge University
Press 1998.
\bibitem{Ko}{\sc H. Kosaki}: {\em Applications of the complex
interpolation method to a von Neumann algebra: Non-commutative
$L_p$-spaces}, J. Funct. Analysis {\bf 56} (1984), 29-78.
\bibitem{La}{\sc E.C. Lance}: {\em Hilbert $C^*$-modules. A toolkit for operator algebraists}, London Math. Soc. Lecture Notes Series 210, Cambridge Univ. Press, Cambridge, 1995.
\bibitem{Mei}{\sc T. Mei}: {\em Operator-valued Hardy spaces},
Preprint.
\bibitem{Mu}{\sc M. Musat}: {\em Interpolation between
non-commutative $BMO$ and non-commutative $L_p$-spaces}, J. Funct.
Analysis {\bf 202} (2003), 195-225.
\bibitem{Ne}{\sc E. Nelson}: {\em Notes on non-commutative integration}, J. Funct. Analysis {\bf 15} (1974), 103-116.
\bibitem{PT}{\sc G. Pedersen, M. Takesaki}: {\em The Radon-Nikodym
theorem for von Neumann algebras}, Acta Math. {\bf 130} (1973),
53-87.
\bibitem{Pet}{\sc K.E. Petersen}: {\em Brownian motion, Hardy
spaces and bounded mean oscillation}, London Math. Soc. Lecture
Notes Series 28, Cambridge Univ. Press, Cambridge 1977.
\bibitem{Pi6}{\sc G. Pisier}: {\em Projections from a von Neumann algebra onto a subalgebra}, Bull. Soc. Math. France{\bf 123} (1995), 139-153.
\bibitem{Pi2}{\sc G. Pisier}: {\em The operator Hilbert space $OH$, complex interpolation and tensor norms}, Mem. Amer. Math. Soc. {\bf 122} (1996), no. 595.
\bibitem{PX}{\sc G. Pisier, Q. Xu}: {\em Non-commutative martingale inequalities}, Comm. Math. Phys. {\bf 189} (1997), no. 3, 667-698.
\bibitem{PX1}{\sc G. Pisier, Q. Xu}: {\em Non-commutative
$L_p$-spaces}, {\em Handbook of the Geometry of Banach Spaces} II,
Ed. W.B. Johnson, J. Lindenstrauss, North Holland (2003),
1459-1517.
\bibitem{Ra}{\sc N. Randrianantoanina}: {\em Non-commutative martingale transforms}, J.
Funct. Analysis {\bf 194} (2002), 181-212.
\bibitem{Ra2}{\sc N. Randrianantoanina}: {\em A weak type inequality for non-commutative martingales and applications}, Preprint.
\bibitem{Ri}{\sc M.A. Rieffel}: {\em Induced representations of
$C^*$-algebras}, Adv. Math. {\bf 13} (1974), 176-257.
\bibitem{Ste1}{\sc E.M. Stein}: {\em Topics in harmonic analysis
related to the Littlewwod-Paley theory}, Princeton University
Press, Princeton, NJ, 1970.
\bibitem {Ste2}{\sc E.M. Stein}: {\em Singular Integrals and
Differentiability Properties of Functions}, Princeton University
Press, Princeton, NJ, 1970.
\bibitem{St}{\sc S. Stratila}: {\em Modular Theory in Operator
Algebras}, Abacus Press, Tunbridge Wells, Kent, England, 1981.
\bibitem{SZ}{\sc S. Stratila, L. Zsido}: {\em Lectures on von
Neumann algebras}, Abacus Press, Tunbridge Wells, Kent, England,
1971.
\bibitem{Ta}{\sc M. Takesaki}: {\em Theory of Operator Algebras I}, Springer-Verlag, New-York, 1979.
\bibitem{Ta1}{\sc M. Takesaki}: {\em Conditional expectations in
von Neumann algebras}, J. Funct. Analysis {\bf 9} (1972), 306-321.
\bibitem{Te}{\sc M. Terp}: {\em $L^p$-spaces associated with von Neumann algebras I and II}, Copenhagen Univ. 1981.
\end{thebibliography}
\end{document}